\documentclass[11pt]{amsart}
\usepackage{multirow,bigdelim}
\usepackage{color}
\usepackage[all,cmtip]{xy}
\usepackage{amsfonts}
\usepackage{dsfont}
\usepackage{amsmath}
\usepackage{extarrows}
\usepackage{mathrsfs}
\usepackage{todonotes}
\usepackage{hyperref}
\usepackage{multimedia}
\usepackage{graphicx}
\usepackage{indentfirst}
\usepackage{bm}
\usepackage{subfigure}
\usepackage{amsthm}
\usepackage{mathrsfs}
\usepackage{latexsym}
\usepackage{amsmath}
\usepackage{amssymb}%»­ÐéÏß¼ýÍ·
\usepackage{microtype}
\usepackage{relsize}
\usepackage{caption}
\usepackage{tikz}
\usepackage{hyperref}
\usepackage{dsfont}%mathds{}Ö§³ÖÊýÑ§×ÖÌå
\usepackage{amssymb}%\mathbb{}Ö§³ÖÊýÑ§×ÖÌå
\usepackage{calc}%ÉèÖÃ¾àÀë£¬ÊµÏÖ¾ÓÖÐ¶Ôxymatrix
\usepackage[marginal]{footmisc}%\section×ó¶ÔÆë
\makeatletter
\def\specialsection{\@startsection{section}{1}%
  \z@{\linespacing\@plus\linespacing}{.5\linespacing}%
  {\normalfont}}% NEW
\def\im{\mathrm{Im}}
\def\section{\@startsection{section}{1}%
  \z@{.7\linespacing\@plus\linespacing}{.5\linespacing}%
  {\normalfont\scshape}}% NEW
\makeatother

\newlength\mylength
\DeclareSymbolFont{SY}{U}{psy}{m}{n}
\DeclareMathSymbol{\emptyset}{\mathord}{SY}{'306}
\newtheorem{theorem}{Theorem}[section]
\newtheorem{corollary}[theorem]{Corollary}
\newtheorem{lemma}[theorem]{Lemma}
\newtheorem{proposition}[theorem]{Proposition}
\theoremstyle{definition}
\newtheorem{definition}[theorem]{Definition}

\theoremstyle{remark}
\newtheorem{remark}[theorem]{Remark}
\theoremstyle{}
\newtheorem{example}[theorem]{Example}

\newcommand{\dr}{\rightrightarrows}
\newcommand{\dl}{\leftleftarrows}
\newcommand{\s}{\simeq}

\newcommand{\pushoutcorner}[1][dr]{\save*!/#1-1.2pc/#1:(-1,1)@^{|-}\restore}
\numberwithin{equation}{section}%\par ÊÇ»»ÐÐÃüÁî
\begin{document}
\captionsetup[figure]{labelfont={bf},labelformat={default},labelsep=space,name={Figure}}
\title{Homotopy Groups and Puppe Sequence of Digraphs}%
\author[Jingyan Li]{Jingyan Li$^{1}$}
%\address{BIMSA}
\author[Jie Wu]{Jie Wu$^{1,3}$}
\author[Shing-Tung Yau]{Shing-Tung Yau$^{1,2}$}
\author[Mengmeng Zhang]{Mengmeng Zhang$^{3,1}$}

\footnote{ $^{1}$ Beijing Institute of Mathematical Sciences and Applications, Beijing 101408, China.}
\footnote{ $^{2}$ Yau Mathematical Sciences Center, Tsinghua
University, Beijing 100084, China.}
\footnote{ $^{3}$ Department of Mathematics, Hebei Normal University,
Shijiazhuang, Hebei 050016, China.}
\footnote{ $E$-$mail\text{ } addresses: \text{ } jingyanli@bimsa.cn(J. Y.\text{ } Li), \text{ } wujie@bimsa.cn(J.\text{ } Wu), \text{ } styau@tsinghua.edu.cn(S.$-$T.\text{ } Yau),\\
 \text{ } mengmengzhang@bimsa.cn(M. M.\text{ } Zhang).$ }
\date{}%
%\footnote{1) Corresponding author}                  %%%   the Fund which you are supported by  %%%
\maketitle%
\begin{abstract}We introduce homotopy groups of digraphs that admit an intuitive description
of grid structures, which is a variation of the GLMY homotopy groups introduced by Grigor'yan, Lin, Muranov and Yau in 2014. This direct approach enables a descriptive interpretation of GLMY theory in applications such as network science. Furthermore, we prove that there exists a long exact
sequence of homotopy groups of digraphs associated to any based digraph map, that is, there exists a digraph version of the Puppe sequence.
\end{abstract}

        % the keywords


% ---------------------------------------------------------------
\section{Introduction}
Detecting and interpreting high-order structures of complex networks is the most challenging problem in
network science \cite{evolu, Beyond, Motifs}.
%which has become a focus topic in some journals such as  Communications Physics.
In terms of mathematics, one possible approach to detect high-order structures
of a complex network given by a (di-)graph is to introduce homology theories of (di-)graphs using techniques
in algebraic topology.
%There have been various mathematical explorations on this topic in topological combinatorics.
One important new homology
theory of digraphs introduced by Grigor'yan-Lin-Muranov-Yau in 2012 \cite{Yau2},
which was originally called path homology and is now called GLMY homology, has achieved various
important applications \cite{Dong, Chow3, Chow2, Shuang}. From the perspective of algebraic topology, the homotopy theory of digraphs plays a more essential role in helping us to understand and gain insight into digraphs.

In 2014, Grigor'yan-Lin-Muranov-Yau introduced the homotopy theory of digraphs \cite{Y-homo}, which coincides with the homotopy theory of graphs in \cite{A-1, A-homo, Balogy, A-2, KL98}, named as $A$-homotopy theory in honor of R. H. Atkin~\cite{Atkin1,Atkin2}. In \cite{Y-homo}, they introduced a notion of $C$-homotopy between loops in digraphs, defined a fundamental group for a based digraph, constructed a loop-digraph $LG$, and recursively defined the $n$-dimensional
homotopy group of a based digraph $G$ by $\pi_n(G):=\pi_1(L^{n-1}G).$  And they demonstrated that 3-cycles and 4-cycles in digraphs may not be contractible when choosing different arrow directions \cite{Y-homo}, compared with the contractability of 3-cycles and 4-cycles in graphs \cite{A-homo}. This observation suggests a certain mystery of digraphs.

Recently, the homotopy group of cubical sets was defined by Carranza and Kapulkin for Kan cubical sets \cite{Carranza} in simplicial approach, which enriching the $A$-homotopy theory of graphs by constructing a functor associating a Kan cubical set to a (simple) graph such that the $A$-homotopy groups of a graph coincide with the homotopy groups of the associated Kan cubical set \cite{cubical}. And by above functor and the loop cubical set $\Omega X$, they obtained the long exact sequence of homotopy groups for fibrations of graphs \cite[Theorem 5.13]{Carranza}.

The purpose of this paper is to lay the foundation for the homotopy theory of digraphs and enrich the GLMY theory \cite{Cover, Yau2, Y-homo, Coho, Y-funda, Steenord, Hoch, Kunne, Yau4, JR}. We introduce the $n$-th homotopy group $\overline{\pi}_n(G)$ for based digraph $G$ by defining the $n$-grid digraph maps, which is a variation of the homotopy group of digraphs introduced in \cite{Y-homo} and the generalization of $A$-homotopy group of graphs in \cite[Proposition 6.5]{A-homo} and cubical sets \cite{cubical, Carranza}. It should be pointed out that our homotopy group is closely connected to the homology group via the Hurewicz theorem, which would
be important for interpreting homology approaches to network science in
detecting the complexity of $n$-dimensional grid structures of data.
%And different from the the loop-digraph $LG$ in \cite{Y-homo} and loop graph $\Omega G$ in \cite{A-1}, we constructed a reduced loop-digraph $\overline{L}G$ with each vertex is a subdivision class of loop in digraph $G$, which satisfies that $\overline{\pi}_n(\overline{L}G)\approx \overline{\pi}_{n+1}(G)$.
Furthermore, different from the simplicial set approach to the fibration theory of (undirected) graphs studied in \cite{cubical, Carranza}, we derive the digraph version of Puppe sequence for arbitrary based digraph map in a canonical constructive approach on the underlying digraphs, which is intimately connected to fibration theory and fibre bundle theory. This direct approach will allow for a descriptive interpretation of GLMY theory in applications such as network science.  An exploration of fibration structures and fibre bundle structures based on underlying digraphs will be one of our subsequent projects. %This result enriches the GLMY theory that has been developed in \cite{Yau2,GLMY,Y-homo,Coho,Y-funda,Steenord,Hoch,Kunne,Yau4,JR}.
  Below we give a detailed description of our main result.

To obtain a grid description of homotopy groups, we modify the definition of homotopy groups of digraphs
by defining the $n$-grid map $$f\colon (I_{m_i},\partial I_{m_i})^{\Box n}\rightarrow (G,\ast)$$ and subdivisions of loops, where $(I_{m_i},\partial I_{m_i})^{\Box n}$ is like the $n$-dimensional cube pair $(I^n,\partial I^n)$ and $m_i$ is the
length of line digraph $I_{m_i}$, $m_i\geq 1$. Then we choose special relative $n$-grid maps $$\overline{f}\colon  (J_{M_i},\partial J_{M_i})^{\Box n}\rightarrow  (G,\ast)$$
 as sphere maps by using subdivisions without changing homotopy type, where $J_M$ is the line digraph of length $M$ with arrow directions alternate. By taking the direct limit $Hom((J,\partial J)^{\Box n};(G,\ast))$ for $\{Hom((J_{M_i},\partial J_{M_i})^{\Box n}; (G,\ast))\}_{M_{i}}$,
the $n$-dimensional homotopy group is defined by $$\overline{\pi}_n(G):= [(J,\partial J)^{\Box n};(G,\ast)].$$
Our main result is the Puppe sequence shown as follows.

\begin{theorem}[Theorem 5.9]
For any based digraph map $f\colon X \rightarrow G$, there is a long exact sequence
\[
\xymatrix@R=0.4cm{
 \cdots \ar[r] & \overline{\pi}_{n+2}(X) \ar[r]^-{f_{n+2}} & \overline{\pi}_{n+2}(G) \ar[r]^-{\Omega\partial_{n+1}} & \overline{\pi}_{n+1}(P_f) \ar[r]^-{\Omega f'_n} & \overline{\pi}_{n+1}(X) \ar[r]^-{\Omega f_n} & \overline{\pi}_{n+1}(G) \\
 \ar[r]^-{\partial_{n+1}} & \overline{\pi}_n(P_f) \ar[r]^-{f'_n} & \overline{\pi}_n(X) \ar[r]^-{f_n} & \overline{\pi}_n(G)
}
\]
of based sets for any $n \geq 0$, where $P_{f}$ is a digraph analogous to the mapping path space in topology. If $n \geq 1$, this is a long exact sequence of groups.
\end{theorem}

The article is organized as follows. In Section 2, we review some basic terminology on the
homotopy theory of digraphs introduced by Yau et al. in \cite{Y-homo}.
 The equivalent description of fundamental group of digraphs in \cite{Y-homo} is given in Section 3 by introducing the notion of subdivision, which inspires us to explore the higher homotopy groups of digraphs. In Section 4, we give an intuitive description of homotopy groups of digraphs by defining subdivision
and prove some properties analogous to classical properties of homotopy groups of topological spaces.
 The Puppe sequence of digraphs will be derived in Section 5 (Theorem \ref{puppe}).

\section{Background}
Based on the homotopy theory introduced by Grigor'yan, Lin, Muranov and Yau in \cite{Y-homo},
we first review some basic definitions and notations in this section.

A \emph{digraph (directed graph)} $G=(V_{G}, A_{G})$ consists of a vertex set $V_{G}$ and an arrow
set $A_{G} \subset\{V_{G} \times V_{G} \backslash \operatorname{diag V_{G}}\}$ of ordered pairs of vertices,
where $\operatorname{diag V_{G}}$ is the diagonal set $\{(v,v)|v\in V_{G}\}$.
 Here, $(v, w) \in A_{G}$ shall be denoted by $v \rightarrow w$.
For convenience we do not distinguish between vertex $v$ and the vertex set $\{v\}$.

Now let us see an important digraph example. Fix $n \geq 0$, denote by $I_n$ a digraph whose
vertex set is $\{0,1,...,n\}$
and arrow set contains exactly one of the arrows $i \rightarrow(i+1)$ or $(i+1) \rightarrow i$
for any $i = 0, 1,...,n-1$. A digraph $I_n$ is called a \emph{line digraph} of length $n$.

%Here another significant definition about digraph is digraph map.
A \emph{digraph map} $f\colon   G\rightarrow H$ is a map $f\colon   V_G \rightarrow V_H $ such that for any
arrow $v\rightarrow w$ in $G$, $f(v)\rightarrow f(w)$ or $f(v)= f(w)$. A digraph $A$ is a \emph{sub-digraph} of a digraph $G$ if $V(A)\subseteq V(G)$ and $E(A)\subseteq E(G)$. A \emph{digraph pair} $(G,A)$ is a digraph $G$ with a
sub-digraph $A$. If $(G,A)$ and $(H,B)$ are digraph pairs,
a \emph{relative digraph map} $f\colon  (G,A)\rightarrow (H,B)$ is a digraph map $f\colon  G\rightarrow H$ such
that $f|_A\colon   A\rightarrow B$ is also a digraph map.

Before we introduce the homotopy of digraph maps,
we need to define the box product of digraphs.
Let $G =(V_G, A_G)$ and $H=(V_H, A_H)$ be two digraphs.
 The \emph{box product} $G\Box H$ is the digraph whose vertex
set is $V_{G}\times V_{H}$ and whose arrow set consists of the arrows $(v,w)\rightarrow (v',w')$ in the cases when $v=v'$ and $w\rightarrow w'$, or $v\rightarrow v'$ and $w=w'$.
%$\{(v,w)\rightarrow (v',w')|\text{either } v=v'\text{ and }w \rightarrow w'\text{ or } v\rightarrow v'\text{ and } w=w'\}$
\subsection{Homotopy of Digraph Maps} Two kinds of homotopy between digraph maps
were introduced by Yau et al. in \cite{Y-homo}. The first homotopy is defined for arbitrary digraph maps, which is generalized by $r$-homotopy in \cite{Asao22, Ser23},
the second one is defined for digraph paths. %Now let us see precise definitions.
\begin{definition}Let $f,~g\colon    (G, A)\rightarrow (H, B)$ be relative digraph maps. We
say that $f$ is \emph{homotopic} to $g$ relative to $A$, denoted by $f\simeq g$\,(rel $A$), if there is a
line digraph $I_n$ and a digraph map $F\colon  G\Box I_n\rightarrow H$ such that $F|_{G\Box\{0\}} = f$,
$F|_{G\Box\{n\}} = g$ and $F|_{A\Box\{i\}}= f|_{A}=g|_{A}$ for any $0\leq i\leq n.$ In particular, if $A=\emptyset,$ we write $f\simeq g.$
%Recently, the $r$-homotopy was introduced by Asao \cite{Asao} and studied in
If $n=1$, we say $f$ is \emph{direct homotopic} to $g$ relative to $A$ and write $f\rightrightarrows g$\,(rel $A$).
\end{definition}
The second homotopy between paths is defined by the mapping cylinder.
Let $G$ be a digraph with base-point $\ast$. A \emph{path} is a relative digraph map
$\phi\colon   (I_n,0)\rightarrow (G,\ast)$. A \emph{loop} on $G$ is a relative digraph map
$\phi\colon  (I_n,\partial I_n)\rightarrow (G,\ast)$, where $\partial I_n$ is the discrete digraph consisting
of the start vertex $0$ and the end vertex $n$ in $I_n$.

For any digraph map $h\colon  G\rightarrow H$, the \emph{cylinder} $C_h$ determined by $h$ is
a digraph whose vertex set is $V_{C_h} = V_G\bigsqcup V_H$ and whose arrow set is
$A_{C_h}= A_G \bigsqcup A_H \bigsqcup \{x\rightarrow h(x)|x\in G\}.$
Similarly there is a \emph{inverse cylinder} $C_h^{-}$ determined by $h$; $C_{h}^{-}$ is the same
as $C_h$ but with $\{x\rightarrow h(x)|x\in G\}$ replaced by $\{h(x)\rightarrow x|x\in G\}$.
From now on, we will use $(0,j)$ and $(1,i)$ to represent the vertices of $G$ and $H$ respectively
in $C_h$ and $C_h^{-}$.

A digraph map $h\colon  I_n\rightarrow I_m$ is called a \emph{shrinking map} if $h(0)=0$, $h(n)=m$ and
$h(i)\leq h(j)$ if $i\leq j$, that is, $h$ is a surjective digraph map preserving vertex order.
\begin{example}\label{shrink}Let $I_3 $ be $\xymatrix@C=0.5cm{ 0\ar[r]& 1& 2 \ar[r]\ar[l]&3}$ and $I_2$ be $\xymatrix@C=0.5cm{0\ar[r]& 1&2.\ar[l]}$ A shrinking map $h\colon  I_3\rightarrow I_2$ is defined by $h(0) = 0$, $h(1)=1$ and $h(2)=h(3) = 2$.
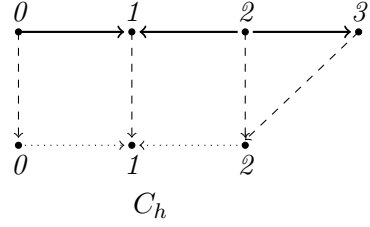
\begin{center}
\begin{tikzpicture}

  % Ìí¼Ó¼ýÍ·
  \draw[->, thick] (1,1.5) -- (2.4,1.5);
   \draw[<-, thick] (2.6,1.5) -- (3.9,1.5);
    \draw[->, thick] (4.1,1.5) -- (5.4,1.5);
\draw[->, dotted] (1,0) -- (2.4,0);
   \draw[<-, dotted] (2.6,0) -- (3.9,0);
   %\draw[->, red, thick] (1,1) -- (1.9,1);
   \draw[->, dashed] (1,1.5) -- (1,0.1);
   \draw[->, dashed] (2.5,1.5) -- (2.5,0.1);
   \draw[->, dashed] (4,1.5) -- (4,0.1);
   \draw[->, dashed] (5.5,1.5) -- (4,0.07);

  % Ìí¼Ó±êÇ©
  \node[below] at (2.75,-0.5) {$C_h$};
  \node[above, black] (2) at (4,1.5) {2};
   \node[above, black] at (1,1.5) {0};
    \node[above, black] at (2.5,1.5) {1};
     \node[above, black] at (5.5,1.5) {3};
     \node[below, black] at (4,0) {2};
   \node[below, black] at (1,0) {0};
    \node[below, black] at (2.5,0) {1};

\filldraw (1,1.5) circle (.04)
(2.5,1.5) circle (.04)
(4 ,1.5) circle (.04)
(5.5,1.5) circle (.04)
(1,0) circle (.04)
(2.5 ,0) circle (.04)
(4,0) circle (.04);
\end{tikzpicture}
\end{center}

%$$\xymatrix@R=1cm{ 0\ar@[red][r]\ar@[green][d]& 1\ar@[green][d]& 2 \ar@[green][d]\ar@[red][l]\ar@[red][r]& 3\ar@[green][dl]\\
%                    0\ar@[blue][r]& 1& 2\ar@[blue][l]& }$$
\end{example}
 \begin{definition}\label{C}Given two paths $\phi\colon   (I_n,0)\rightarrow (G,\ast)$ and $\psi\colon   (I_m,0)\rightarrow (G,\ast)$,
we say $\phi$ is \emph{one-step direct $C$-homotopic} to $\psi$, denoted by $\phi\simeq^1\psi$,
if there exists a shrinking map $h\colon   I_n\rightarrow I_m$ and a digraph map $F\colon   C_h\rightarrow G$ such
that $F|_{I_n}= \phi$ and $F|_{I_m}=\psi$. \\
We say $\phi$ is \emph{one-step inverse $C$-homotopic} to $\psi$, denoted by $\phi\simeq^{-1}\psi$,
 if there exists a shrinking map $h\colon   I_n\rightarrow I_m $ and a digraph map $F\colon  C_h^{-}\rightarrow G$ such
that $F|_{I_n} = \phi$ and $F|_{I_m} = \psi$.
\end{definition}
\begin{example}
Let $h$ be the shrinking map in Example \ref{shrink}, $\phi\colon  (I_3,\partial I_3)\rightarrow (G,\ast)$ is the thickened digraph map and $\psi\colon  (I_2,\partial I_2)\rightarrow (G,\ast)$ is the dotted digraph map, as illustrated below. Then $F\colon  C_h\rightarrow G$ is the digraph map sending the arrow to the arrow with same shape or vertex.
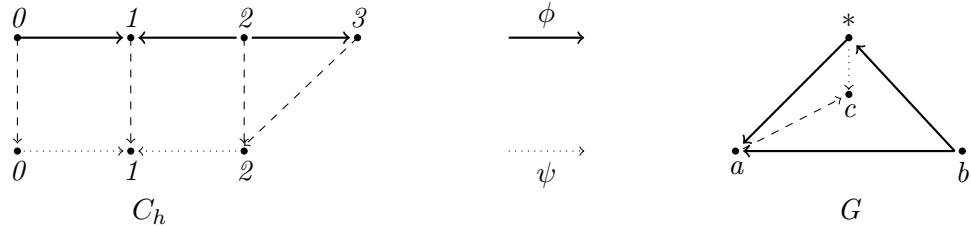
\begin{center}
\begin{tikzpicture}

  % Ìí¼Ó¼ýÍ·
  \draw[->, thick] (13.4,0) -- (12.1,1.4);
  \draw[->, thick] (13.4,0) -- (10.6,0);
   \draw[->, thick] (12,1.5) -- (10.6,0.1);
  \draw[->, thick] (1,1.5) -- (2.4,1.5);
   \draw[<-, thick] (2.6,1.5) -- (3.9,1.5);
    \draw[->, thick] (4.1,1.5) -- (5.4,1.5);
\draw[->, dotted] (1,0) -- (2.4,0);
   \draw[<-, dotted] (2.6,0) -- (3.9,0);
   %\draw[->, red, thick] (1,1) -- (1.9,1);
   \draw[->, dashed] (1,1.5) -- (1,0.1);
   \draw[->, dashed] (2.5,1.5) -- (2.5,0.1);
   \draw[->, dashed] (4,1.5) -- (4,0.1);
   \draw[->, dashed] (5.5,1.5) -- (4,0.07);
   \draw[->, thick] (7.5,1.5) -- (8.5,1.5);
\draw[->, dotted] (7.5,0) -- (8.5,0);
\draw[<-, dotted] (12,0.81) -- (12,1.4);
\draw[->, dashed]  (10.6,0.03) -- (11.9,0.7);

  % Ìí¼Ó±êÇ©
  \node[below] at (2.75,-0.5) {$C_h$};
  \node[above, black] (2) at (4,1.5) {2};
   \node[above, black] at (1,1.5) {0};
    \node[above, black] at (2.5,1.5) {1};
     \node[above, black] at (5.5,1.5) {3};
     \node[below, black] at (4,0) {2};
   \node[below, black] at (1,0) {0};
    \node[below, black] at (2.5,0) {1};
    \node[above, black] at (8,1.5) {$\phi$};
    \node[below, black] at (8,0) {$\psi$};
    \node[below, black] at (10.5,0) {a};
    \node[below, black] at (13.5,0) {b};
    \node[above, black] at (12,1.5) {$\ast$};
    \node[below, black] at (12,0.75) {c};
     \node[below, black] at (12,-0.5) {G};

\filldraw (1,1.5) circle (.04)
(2.5,1.5) circle (.04)
(4 ,1.5) circle (.04)
(5.5,1.5) circle (.04)
(1,0) circle (.04)
(2.5 ,0) circle (.04)
(4,0) circle (.04)
(10.5,0) circle (.04)
(13.5,0) circle (.04)
(12,1.5) circle (.04)
(12,0.75) circle (.04);
\end{tikzpicture}
\end{center}
%\begin{center}
%$\xymatrix@R=1.0cm{ 0\ar@[red][r]\ar@[green][d]& 1\ar@[green][d]& 2 \ar@[green][d]\ar@[red][l]\ar@[red][r]& 3\ar@[green][dl]&  &   \textcolor[rgb]{1.00,0.00,0.00}{\phi}\ar@[red][r] & & \\
%                    0\ar@[blue][r]& 1& 2\ar@[blue][l]& &  &  \textcolor[rgb]{0.10,0.05,0.95}{\psi}\ar@[blue][r]& & }$
%                     $\xymatrix@R=0.4cm{ & & \ast\ar@[red][ddll] \ar@[blue][d] & &  \\
%                    &&c\ar@[green][dll] \ar@[blue][u] & &  \\
%                     a& & &  &b\ar@[red][llll]^G\ar@[red][uull] }$
%\end{center}
\end{example}
To facilitate the definition of C-homotopy given below, we first introduce a piece of notation: we write \( f \simeq^1 g \) or \( g \simeq^{-1} f \) as \( f \rightarrow g \).
\begin{definition}
Let $\phi\colon  (I_n,0)\rightarrow (G,\ast)$ and $\psi\colon  (I_m,0)\rightarrow (G,\ast)$ be two paths.
We say $\phi$ is \emph{$C$-homotopic} to $\psi$, denoted by $\phi\simeq^C \psi,$ if there
exists a finite sequence of paths $\{\phi_j\}_{j=0}^{l}$ such that $\phi_0 =\phi$,
$\phi_l = \psi$ and $\phi_j \rightarrow \phi_{j+1}$ or $\phi_j\leftarrow \phi_{j+1}$ for $j=0,\cdots,l-1$.
\end{definition}
Let $\phi$ and $\psi$ be two loops, $\phi$ is \emph{$C$-homotopic} to $\psi$ if $\phi$ is C-homotopic
 to $\psi$ as paths and $\{\phi_j\}_{j=0}^{l}$ is a sequence of loops. One can easily verify that $C$-homotopy of paths and loops gives an equivalence relation. Following this, Yau and his coauthors defined the fundamental group of $G$.

Assume that $\phi\colon  (I_n,0)\rightarrow (G,\ast)$ and $\psi\colon (I_m,0)\rightarrow (G,\ast)$ are two paths. The \emph{concatenation} of $\phi$ and $\psi$ is the path $ \phi\vee \psi\colon (I_{m+n},0)\rightarrow (G,\ast)$ such that $\phi\vee \psi|_{[0,n]}=\phi$ and $\phi\vee \psi|_{[n,m+n]}=\psi$. In particular, the concatenation of two loops is still a loop. Moreover, the \emph{inverse} of loop $\gamma\colon (I_l,\partial I_l)\rightarrow (G,\ast)$ is the loop $\gamma^{-1}\colon (\widehat{I}_{l},\partial \widehat{I}_l)\rightarrow (G,\ast)$ such that $\gamma^{-1}(i)=\gamma(l-i)$, where $\widehat{I}_{l}$ is the line digraph of length $l$ with arrow $i\rightarrow j$ if and only if $l-i \rightarrow l-j$ is in $I_l$.

\begin{definition}
  Let $G$ be a digraph with base-point $\ast$. The \emph{fundamental group} $\pi_1(G)$ is the
group consisting of $ C$-homotopy equivalence classes of loops in $G$ with multiplication being
concatenation of loops.
\end{definition}
    \subsection{GLMY Higher Homotopy Group \texorpdfstring{$\pi_n(G)$}{ } } Let $G$ be a digraph with base-point $\ast$ and $V_2$ be a discrete digraph having two vertices $0$ and $1$. Then we shall write $Hom((V_2,0);(G,\ast))=\{f\colon (V_2,0)\rightarrow (G,\ast)\}.$

Recall that $f$ is \emph{homotopic} to $g$ if
there exists a line digraph $I_n$ and a digraph map $F\colon  I_n\rightarrow G$ such that $F(0)=f(1)$ and $F(n)=g(1)$ for any two maps $f$ and $g$ in $Hom((V_2,0);(G,\ast))$.

Define $\pi_0(G)$ as the set of homotopy classes of $Hom((V_2,0);(G,\ast))$.
It is not difficult to see that $Hom((V_2,0);(G,\ast))$ is isomorphic to $Hom(0,G)=\{\phi\colon  0\rightarrow G\}$. Obviously $\pi_0(G)$ is exactly the set of homotopy classes
of $Hom(0,G)$. Then $\pi_0(G)$ coincides with the set of path-components of $G$.

Yau et al. also defined a loop-digraph $LG$ such that $\pi_1(G)=\pi_0(LG)$ for any based digraph $G$,
just like the loop space of a topological space. The \emph{loop-digraph} $LG$ is a
digraph whose vertex set consists of all loops in $G$ and whose arrow set contains the
arrow $f\rightarrow g$ if and only if $f\simeq^1 g$ or $g\simeq^{-1}f$. By
definition of $LG$, it is easy to see that $\pi_0(LG)\approx \pi_1(G).$

Based on $\pi_1(G)=\pi_0(LG)$, Yau et al.
defined the higher homotopy group inductively.
\begin{definition}
Let $G$ be a digraph with base-point $\ast$. The \emph{$n$-dimensional homotopy group $\pi_n(G)$} is defined by $$\pi_n(G):=\pi_{n-1}(LG)$$ for $n\geq1$.
\end{definition}

For convenience, we will refer digraph $G$ and digraph map $f $ as the
 based digraph and based digraph map and $I_n$ as the digraph with base-point $0$ if there is no other specification. That is to say, the digraph map $f\colon  I_n\rightarrow G$ always is a path.
\section{Equivalent Description of the Fundamental Group}
The purpose of this section is to give an equivalent description for $C$-homotopy between loops and the fundamental group of a based digraph.
By using this description of $C$-homotopy, we will proceed to define the higher homotopy groups of a digraph. To begin, we recall some basic definitions.
%Moreover, we will also give a simpler representation of the fundamental group of a based digraph.

A \emph{directed set}~\cite{Mac} is a poset $(A,\leq)$ such that for
any $\alpha$, $\beta\in A$, there is an upper bound $\gamma\in A$:  $\alpha \leq \gamma$, $\beta\leq \gamma.$ Let $\mathcal{C}$ be a category. A \emph{directed system} $\{J_{\alpha}; j_{\alpha}^{\beta}\}$ is a family of objects $\{J_{\alpha}| J_{\alpha}\in \mathcal{C}\}_{\alpha\in A}$ and morphisms $j_{\alpha}^{\beta}\colon  J_{\alpha}\rightarrow J_{\beta} $ for all $\alpha \leq \beta$ such that $j_{\beta}^{\gamma}\circ j_{\alpha}^{\beta}= j_{\alpha}^{\gamma}$ if $\alpha\leq \beta\leq \gamma$ and $j_{\alpha}^{\alpha}=1$, where $A$ is a directed set.

Consider a category $\mathcal{C}$ of sets. It is well known that direct limits (i.e. filtered colimits) exist in $\mathcal{C}$. Given a directed system $\{J_{\alpha}; j_{\alpha}^{\beta}\}$ indexed by a directed set $(A,\leq)$, the \emph{direct limit} 
$\lim\limits_{\rightarrow}J_{\alpha}$ is explicitly given by  $$\bigsqcup\limits_{\alpha\in A}J_{\alpha}/\sim,$$ 
where $\bigsqcup$ denotes the disjoint union of the sets $J_{\alpha}$, and the equivalence relation `$\sim$' is defined as follows: for $x_{\alpha}\in J_{\alpha}$ and $x_{\beta}\in J_{\beta}$, we say that 
$x_{\alpha}\sim x_{\beta}$ if and only if there exists a $\gamma\in A$ with $\alpha\leq \gamma$ and
$\beta\leq \gamma$ such that $$j_{\alpha}^{\gamma}(x_{\alpha}) = j_{\beta}^{\gamma}(x_{\beta}).$$  Furthermore, there is a family of canonical injections $$i_{\alpha}\colon  J_{\alpha} \rightarrow \lim\limits_{\rightarrow}J_{\alpha},\quad x\mapsto \{ x\}$$
for each $\alpha\in A$ such that $i_{\beta}\circ j_{\alpha}^{\beta}=i_{\alpha}.$ 

Now we give a kind of special line digraph $J_n$. For each non-negative integer $n$, the digraph $J_n$ is the line digraph of length $n$, with arrows alternating as $i \to i+1$ if $i$ is even, and $i \leftarrow i+1$ if $i$ is odd, for all $0 \leq i \leq n-1$. For example, 
\begin{center}
\begin{tikzpicture}
\node[below, black] at (7.4,0.27) {$J_3:$};
\node[below, black] at (11,0) {3.};
   %\node[above, black] at (6.5,0.2) {h};
    \node[below, black] at (8,0) {0};
     \node[below, black] at (9,0) {1};
     \node[below, black] at (10,0) {2};
     
\filldraw (8,0) circle (.04)
(9,0) circle (.04);
\filldraw (11,0) circle (.04);
\filldraw (10,0) circle (.04);
   \draw[->, thick] (8.1,0) -- (8.9,0);
   \draw[->, thick] (9.9,0) -- (9.1,0);
   \draw[->, thick] (10.1,0) -- (10.9,0); 
\end{tikzpicture}
\end{center}
In contrast to the line digraph $I_n$, the arrows of $J_n$ are completely determined by its length $n$. In what follows, we will refer to $J_n$ as the standard line digraph, and to the paths (loops) from $J_n$ as the standard paths (loops).

Observe that the natural number $(\mathbb{N},\leq)$ forms a directed set. Then indexed by $\mathbb{N}$, the standard line digraphs $\{J_m\}_{m\geq 0}$ together with the canonical inclusion maps $j_m^n\colon J_m \to J_n$  form a directed
system $\{J_m; j_{m}^{n}\}_{m\leq n}$. The direct limit $\lim\limits_{\rightarrow}J_{m}$ for this system, denoted by $J_{\infty}$, is
$\bigcup\limits_{m\geq 0}J_{m}$, where $\bigcup$ is the union of the underlying vertex sets and arrow sets of the $J_m$.  In other words, \(J_\infty\) is the infinite standard line digraph obtained by gluing together the finite \(J_m\) along their canonical inclusions.  
By using the standard line digraphs, we can prove that every path admits a standard representative within its $C$-homotopy class.

\begin{lemma}\label{standard}
Let $G$ be a based digraph. For any path $f\colon I_m\rightarrow G$, there is a standard path $\widehat{f}\colon  J_M\rightarrow G$ for some $M\geq m$ such that $\widehat{f}\simeq^1 f$.
\begin{proof}

To prove this lemma, we construct three digraph maps: $\widehat{f}$, a shrinking map $h\colon J_{M} \to I_{m}$, and $F\colon C_{h} \to G$, where $C_{h}$ denotes the cylinder of $h$. The case $m=0$ is immediate, so assume $m\geq 1$.

Let $\widehat{f}(0) = f(0)$ and $h(0) = 0$, and extend $\widehat{f}$ and $h$ inductively by comparing the $i$-th arrow in $I_m$ with the $j$-th arrow in $J_{\infty}$, starting at $i=j=1$: \begin{itemize} \item If the arrows agree in direction, let $\widehat{f}(j) = f(i)$ and $h(j) = i$, then increment both $i$ and $j$ by 1. \item If the arrows differ, let $\widehat{f}(j) = f(i-1)$, $\widehat{f}(j+1) = f(i)$, $h(j) = i-1$, $h(j+1) = i$, then increment $i$ by 1 and $j$ by 2. \end{itemize} Repeat this process until reaching the last arrow of $I_m$. Let $M=j$, and set $\widehat{f}(M) = f(m)$ and $h(M) = m$. Then $\widehat{f}$ restricts to a digraph map $J_{M} \to G$.

By construction, $h$ is a shrinking map. Finally, we define the digraph map $F\colon C_{h} \to G$ by extending $\widehat{f}$ and $f$ along the cylinder: on each arrow $(0,j) \to (1,i)$ with $h(j) = i$, set $F\left( (0,j) \to (1,i) \right) = \widehat{f}(j) = f(i)$. Then $F$ restricts to $\widehat{f}$ on $0\Box J_M$ and to $f$ on $1\Box I_m$, establishing a $1$-homotopy between $\widehat{f}$ and $f$.
\end{proof}
\end{lemma}
According to this lemma, there is an interesting result about the fundamental group of a digraph.
First let us make some preliminary statements.
Denote the set $\{f\colon  (J_m,\partial J_m)\rightarrow (G,\ast)\}$ by $Hom((J_m,\partial J_m);(G,\ast))$.

Based on
the directed set $\mathbb{N}$, there is a directed
system $\{Hom((J_m,\partial J_m);(G,\ast)); $ $l_{m}^{n}\}_{m\leq n}$,
where $$l_{m}^{n}\colon Hom((J_m,\partial J_m);(G,\ast))\rightarrow Hom((J_n,\partial J_n);(G,\ast)),\text{ }f\mapsto \overline{f},$$
with $ \overline{f}(i)=\left\{
                     \begin{array}{ll}
                       f(i), & \hbox{$i\leq m$,} \\
                       f(m), & \hbox{$i> m$.}
                     \end{array}
                   \right
.$

Since $Hom((J_m,\partial J_m);(G,\ast))$ is an object of the category $\mathcal{S}et$, then
$$\lim\limits_{\rightarrow}Hom((J_m,\partial J_m);(G,\ast))= \bigsqcup\limits_{m\geq0}Hom((J_m,\partial J_m);(G,\ast))/\sim,$$
where $f_{m}\sim f_{n}$ if and only if there exists an integer $k\in \mathbb{N}$ with $m\leq k$ and $n\leq k$
such that $l_{m}^{k}(f_{m}) = l_{n}^{k}(f_{n})$ for any $f_{m}\in Hom((J_m,\partial J_m);(G,\ast))$ and
$f_{n}\in Hom((J_n,\partial J_n);$ $(G,\ast))$. From the definition of direct limit, there is a family of injections 
$$\{i_m\colon \operatorname{Hom}((J_m,\partial J_m);(G,\ast))\rightarrow \lim\limits_{\rightarrow}\operatorname{Hom}((J_m,\partial J_m);(G,\ast)) \}$$ such that $i_n\circ l_{m}^{n} = i_m.$

In the proof of Lemma \ref{standard}, the construction of the shrinking map
$h\colon  J_{M}\rightarrow I_{n}$ implies that there is a map $$\Gamma_{I_m}\colon \operatorname{Hom}((I_m,\partial I_m);(G,\ast))\rightarrow \bigcup \limits_{M \geq 1}\operatorname{Hom}((J_M,\partial J_M);(G,\ast)) \stackrel{i_M}{\longrightarrow} \lim\limits_{\rightarrow}\operatorname{Hom}((J_M,\partial J_M);(G,\ast))$$ sending
 $f$ to $i_M(\widehat{f})= \{\widehat{f}\}$ for any $I_m$. Furthermore, the collection of the maps $\Gamma_{I_m}$ determines a map $$\Gamma\colon  \bigcup \limits_{I_m\in \mathcal{I}}\operatorname{Hom}((I_m,\partial I_m);(G,\ast))\rightarrow \lim\limits_{\rightarrow}\operatorname{Hom}((J_M,\partial J_M);(G,\ast)),$$ where $\mathcal{I}$ is the set of all line digraphs.

Now let us define the $C$-homotopy in $\lim\limits_{\rightarrow}\operatorname{Hom}((J_M,\partial J_M);(G,\ast)).$ We say $\{\widehat{f}\} \s^C \{\widehat{g}\}$ if and only if $f\s^C g$. One can easily check that this is well-defined regardless of the representatives in the same direct limit class. 
% Assume that since $\widehat{f}\s^1 \widehat{f}_{min}$, $\widehat{g}\s^1 \widehat{g}_{min}$ and $\widehat{f}_{min}\s^C \widehat{g}_{min}$, then $\widehat{f}\s^C \widehat{g}.$ That is to say, $\{\widehat{f}\} \s^C \{\widehat{g}\}$ if and only if $\widehat{f}\s^C \widehat{g}.$
In this way we obtain a set $[J,G]^{C}$ consisting of the $C$-homotopy classes of the maps in $\lim\limits_{\rightarrow}\operatorname{Hom}((J_M,\partial J_M);(G,\ast))$. Then the $C$-homotopy class of $\{f\}$ is denoted by $[\{f\}]$.

Moreover, $[J,G]^{C}$ forms a group with the multiplication given by the concatenation of loops.
To guarantee that the concatenation of loops is still a standard loop,  for any $\{f\}$, we choose
 $f\colon (J_{m},\partial J_m)\rightarrow (G,\ast)$ with even length as the representatives of $\{f\}$. The  multiplication $\mu$ is defined as follows $$\mu\colon [J,G]^{C}\times [J,G]^{C}\rightarrow [J,G]^{C}, \text{ }([\{f\}],[\{g\}])\mapsto [\{f\vee g\}] $$ on $[J,G]^{C}$, where $f\colon (J_{m},\partial J_m)\rightarrow (G,\ast)$ and $g\colon (J_{n},\partial J_n)\rightarrow (G,\ast)$ are standard loops of even length. We claim that the concatenation of $[\{f\}]$ and $[\{g\}]$ is independent of the choice of the representation elements $f$ and $\{f\}$. Suppose $\{f\} = \{f'\}$ and $f$ is longer than $f'$, then $f\vee g \s^1 f'\vee g$, and if $f$ is shorter than $f'$, then $f'\vee g \s^{1} f\vee g$. Thus $f\vee g \s^C f'\vee g.$ Also, if $\{f\}\s^C \{f'\}$, that is,  $f\s^C f'$, then $f\vee g \s^C f'\vee g$, that is, $[\{f\vee g\}] = [\{f' \vee g\}].$ Hence $\mu$ is well-defined and associative.

Obviously, the inverse loop $f^{-1}$ of the standard loop $f\colon (J_m,\partial J_m)\rightarrow (G,\ast)$ is also a standard loop with even length, where $$f^{-1}\colon  (J_m,\partial J_m)\rightarrow (G,\ast), \text{ }i\mapsto f(m-i).$$

The unit element is $1\colon J_{2}\rightarrow \ast$ in $[J,G]^{C}$. Thus $[J,G]^{C}$ is a group. In what follows, we will identify $\{\widehat{f}\}$ with $\widehat{f}$ and identify the $C$-homotopy class $[\{\widehat{f}\}]$
of $\{\widehat{f}\}$ with $[\widehat{f}]$. We write $\mu([f],[g])$ as $ [f]\cdot [g]$.

\begin{theorem}\label{funda}
Let $G$ be a digraph with base-point $\ast$. The map $\Gamma$ induces an isomorphism
$$\Gamma_1\colon \pi_1(G)\rightarrow [J,G]^{C},\text{ }[f]\mapsto [\widehat{f}].$$
\begin{proof}
Recall the definition of $\Gamma$, for any loop $f\colon (I_{m},\partial I_{m})\rightarrow (G,\ast)$,
$\Gamma(f) = \{\widehat{f}\}$. If $[f_1] = [f_{2}]$ in $\pi_1(G)$, that is, $f_1\s^{C}f_2$, then we have
$\widehat{f}_1 \simeq^{1}f_1 \simeq^{C} f_{2}$ and $\widehat{f}_{2} \simeq^{1}f_{2}$.
Hence $\widehat{f}_1 \simeq^{C} \widehat{f}_{2}$, that is, $\{\widehat{f}_1\}\s^{C} \{\widehat{f}_2\}.$
Thus $\Gamma_{1}$ is well-defined.

For any elements $[f_1],\text{ }[f_2]$ in $ \pi_1(G)$, if $\Gamma(f_{1})\simeq^{C} \Gamma(f_{2})$,
then $\widehat{f}_1\s^C \widehat{f}_2.$ By $\widehat{f}_{1} \simeq^{1} f_{1}$ and $\widehat{f}_{1}\simeq^{C} \widehat{f}_{2}\simeq^{1} f_{2}$,
we have $f_1 \simeq^C f_2$. So $\Gamma_{1}$ is injective. Clearly, $\Gamma_{1}$ is surjective.

 Next we check that $\Gamma_{1}$ is a homomorphism. For any $C$-homotopy classes $[f]$
and $[g]$ in $\pi_{1}(G)$, we have
$\Gamma_{1}([f]\cdot[g])= \Gamma_{1}([f\vee g])= [\widehat{f\vee g}] = [f\vee g] = [f]\cdot[g]=[\widehat{f}]\cdot[\widehat{g}]= \Gamma_{1}([f])\cdot \Gamma_{1}([g]).$
Thus $\Gamma_1$ is an isomorphism.
\end{proof}
\end{theorem}
This theorem states that if you want to study the fundamental group of a based digraph, you just need to think about all the standard loops, not all loops.

Next we introduce a new equivalent condition for $C$-homotopy, which will
inspire the definition of the higher homotopy groups. This requires a significant definition.
\begin{definition}
Let $f\colon  I_m\rightarrow G$ be a path and let
$h\colon I_M\rightarrow I_m$ be a shrinking map. The digraph map $ \overline{f}= f\circ h\colon I_M\rightarrow G$ is called
a \emph{subdivision} of $f$.
\end{definition}
The notion of subdivision is the digraph version of shomotopy between loops in simplicial complex introduced by R. H. Atkin in \cite{Atkin2} and extension of loops in graph in \cite{Y-homo} with considering the arrow direction in line digraph. And it should be highlighted that we will construct a reduced loop-digraph $\overline{L}G$ later by using this subdivision. Now let us see the following example.
\begin{example}
Let $f\colon  (I_3,\partial I_3) \rightarrow (G,\ast)$. Given a shrinking map \( h\colon I_5 \rightarrow I_3 \) that sends the dotted arrows \( 1 \leftarrow 2 \) and \( 4 \to 5 \) to the vertices \( 1 \) and \( 3 \), respectively, as shown in the following diagram.

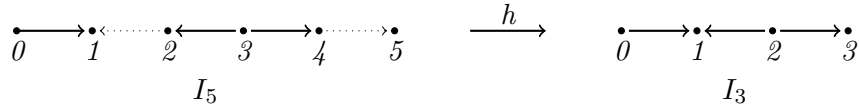
\begin{center}
\begin{tikzpicture}
  % Ìí¼Ó¼ýÍ·
\draw[->, thick] (0,0) -- (0.9,0);
   \draw[<-, dotted] (1.1,0) -- (1.9,0);
   %\draw[->, red, thick] (1,1) -- (1.9,1);
   \draw[->, thick] (2.9,0) -- (2.1,0);
   \draw[->, thick] (3.1,0) -- (3.9,0);
   \draw[->, dotted] (4.1,0) -- (4.9,0);
   \draw[->, thick] (6,0) -- (7,0);
   \draw[->, thick] (8.1,0) -- (8.9,0);
   \draw[->, thick] (9.9,0) -- (9.1,0);
   \draw[->, thick] (10.1,0) -- (10.9,0);

  % Ìí¼Ó±êÇ©
  \node[below] at (2.5,-0.5) {$I_5$};
  \node[below, black] at (2,0) {2};
   \node[below, black] at (6.5,0.5) {h};
    \node[below, black] at (5,0) {5};
     \node[below, black] at (3,0) {3};
     \node[below, black] at (4,0) {4};
   \node[below, black] at (1,0) {1};
    \node[below, black] at (0,0) {0};
    \node[below] at (9.5,-0.5) {$I_3$};
  \node[below, black] at (11,0) {3};
   %\node[above, black] at (6.5,0.2) {h};
    \node[below, black] at (8,0) {0};
     \node[below, black] at (9,0) {1};
     \node[below, black] at (10,0) {2};

\filldraw (0,0) circle (.04)
(8,0) circle (.04);
\filldraw (2,0) circle (.04)
(9,0) circle (.04)
(1,0) circle (.04);
\filldraw 
(4,0) circle (.04)
(5 ,0) circle (.04)
(11,0) circle (.04);
\filldraw (3,0) circle (.04)
(10,0) circle (.04);
\end{tikzpicture}
\end{center}
%$$ \xymatrix@R=0.2cm@C=0.5cm{\textcolor[rgb]{1.00,0.00,0.00}{ 0} \ar[r] & \textcolor[rgb]{0.50,1.00,0.00}{1} \ar@[green][r]& \textcolor[rgb]{0.50,1.00,0.00}{2} & \textcolor[rgb]{0.00,0.50,1.00}{3} \ar[l]\ar[r] & \textcolor[rgb]{1.00,0.00,0.00}{4}\ar@[red][r]& \textcolor[rgb]{1.00,0.00,0.00}{5} \\
%& &I_5 & &}\stackrel{h}{\longrightarrow}\xymatrix@R=0.2cm@C=0.5cm{
%  \textcolor[rgb]{1.00,0.00,0.00}{0} \ar[r] & \textcolor[rgb]{0.50,1.00,0.00}{1} & \textcolor[rgb]{0.00,0.50,1.00}{2} \ar[l] \ar[r] & \textcolor[rgb]{1.00,0.00,0.00}{3},\\
%  &I_3& }$$
Then
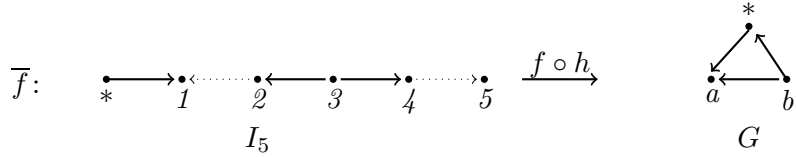
\begin{center}
\begin{tikzpicture}
  % Ìí¼Ó¼ýÍ·
\draw[->, thick] (1,0) -- (1.9,0);
   \draw[<-, dotted] (2.1,0) -- (2.9,0);
   %\draw[->, red, thick] (1,1) -- (1.9,1);
   \draw[->, thick] (3.9,0) -- (3.1,0);
   \draw[->, thick] (4.1,0) -- (4.9,0);
   \draw[->, dotted] (5.1,0) -- (5.9,0);
   \draw[->, thick] (6.5,0) -- (7.5,0);
   \draw[->, thick] (9.5,0.65) -- (9,0.1);
   \draw[->, thick] (9.9,0) -- (9.1,0);
   \draw[->, thick] (10,0) -- (9.6,0.6);

  % Ìí¼Ó±êÇ©
  %\node[below] at (2.5,-0.5) {$I_5$};
  \node[below, black] at (0,0.3) {$\overline{f}\colon$};
  \node[below, black] at (7,0.5) {$f\circ h$};
  \node[below, black] at (3,0) {2};
    \node[below, black] at (6,0) {5};
     \node[below, black] at (4,0) {3};
     \node[below, black] at (5,0) {4};
   \node[below, black] at (2,0) {1};
    \node[below, black] at (1,0) {$\ast$};
    \node[below] at (9.5,-0.5) {$G$};
    \node[below] at (3,-0.5) {$I_5$};
  %\node[below, black] at (11,0) {$\ast .$};
   %\node[above, black] at (6.5,0.2) {h};
    \node[above, black] at (9.5,0.7) {$\ast$};
     \node[below, black] at (9,0) {a};
     \node[below, black] at (10,0) {b};
\filldraw (1,0) circle (.04)
(9.5,0.7) circle (.04);
\filldraw (3,0) circle (.04)
(9,0) circle (.04)
(2,0) circle (.04);
\filldraw 
(5,0) circle (.04)
(6 ,0) circle (.04);
\filldraw (4,0) circle (.04)
(10,0) circle (.04);
\end{tikzpicture}
\end{center}
%$$\overline{f}= f \circ h:\xymatrix@C=0.5cm@R=0.2cm{
% \textcolor[rgb]{1.00,0.00,0.00}{ 0} \ar[r] & \textcolor[rgb]{0.50,1.00,0.00}{1} \ar@[green][r]& \textcolor[rgb]{0.50,1.00,0.00}{2} & \textcolor[rgb]{0.00,0.50,1.00}{3} \ar[l]\ar[r] & \textcolor[rgb]{1.00,0.00,0.00}{4}\ar@[red][r]& \textcolor[rgb]{1.00,0.00,0.00}{5}\\
% &&I_5&&}\stackrel{h}{\longrightarrow}\xymatrix@C=0.5cm@R=0.2cm{
%  \textcolor[rgb]{1.00,0.00,0.00}{0} \ar[r] & \textcolor[rgb]{0.50,1.00,0.00}{1} & \textcolor[rgb]{0.00,0.50,1.00}{2} \ar[l] \ar[r] & \textcolor[rgb]{1.00,0.00,0.00}{3} \\
%  &I_3&}\stackrel{f}{\longrightarrow}\xymatrix@C=0.5cm@R=0.2cm{
%  \textcolor[rgb]{1.00,0.00,0.00}{\ast} \ar[r] & \textcolor[rgb]{0.50,1.00,0.00}{a} & \textcolor[rgb]{0.00,0.50,1.00}{b} \ar[l] \ar[r] & \textcolor[rgb]{1.00,0.00,0.00}{\ast} .\\
%  &G&}$$
\end{example}
An important property of subdivisions is the following.
\begin{proposition}\label{subd}
   For any path $f\colon I_m\rightarrow G$ and any subdivision $\overline{f}$ of $f$, we have $\overline{f}\s^1 f$.
\begin{proof}
Suppose that $\overline{f}$ is a subdivision of $f$ by a shrinking map $h\colon I_M\rightarrow I_m$, that is,
$\overline{f} = f\circ h$. We construct a map $F\colon C_h\rightarrow G$ such that
$F|_{\{0\}\Box I_M}= \overline{f}$ and $F|_{\{1\}\Box I_m}= f$. Observe that $(0,j)\rightarrow (1,i)$ is an arrow of
$C_h$ if and only if $h(j)=i$. For such an arrow, define $F((0,j)\rightarrow (1,i)) = f(i)$. Clearly, $F$ is a digraph map.
Hence $\overline{f}\s^{1} f.$
\end{proof}
\end{proposition}
Subdivision lets us see $C$-homotopy in a new light, which can then be applied to define
higher homotopy groups. To give the precise description, we introduce the following notation.

\begin{definition}\label{1Fhomo}
Let $f\colon I_{m}\rightarrow G$ and $g\colon I_{n}\rightarrow G$ be two paths in $G$. We call $f $ \emph{one-step $F$-homotopic} to $g$ if there exist subdivisions $\overline{f}$ and $\overline{g}$ of $f$ and $g$ respectively such that $f\dr g$, denoted by $f\simeq_{1}g$. We also call $g$ \emph{one-step inverse $F$-homotopic} to $f$ and write $g\s_{-1}f$.

 More generally, we say that $f$ is \emph{$F$-homotopic} to $g$, denoted by $f\simeq_F g$  if there is a finite sequence of $\{f_i\}_{i=0}^l$ such that $f_{0}=f$, $f_{l}=g$ and there are one-step $F$-homotopies $f_i\s_{1}f_{i+1}$ or $f_i\s_{-1}f_{i+1}$ for $0\leq i\leq l-1$.
\end{definition}

Obviouly, $F$-homotopy is an equivalence relation between paths, and we will show that in Lemma \ref{equ} it is equivalent to the $C$-homotopy equivalence relation in Definition \ref{C}. 

\begin{lemma}\label{equ}
 For any two standard paths $f\colon  J_m\rightarrow G$ and $g\colon J_n \rightarrow G$, $f\simeq^C g$ if and only if $f\simeq_F g$.
 \begin{proof}
Let $f\colon  J_m\rightarrow G$ and $g\colon  J_n\rightarrow G$ be paths. We only need to show that if $f\s_1 g$ or $f\s_{-1} g$, then $f\s^C g$. For the case $f\s_1 g$, there exist subdivisions $\overline{f}$ and $\overline{g}$ of $f$ and $g$ respectively such that $\overline{f} \dr \overline{g}$, so we have $\overline{f} \s^{1} \overline{g}\s^1 g$ and $\overline{f}\simeq^{1} f$, and therefore $f\s^C g$. Similarly if $f\s_{-1} g$, we also have $f\s^C g$.

Conversely, suppose that $f\s^C g$. It is sufficient to consider two cases: $f\rightarrow g$ and $f\leftarrow g$, which can be reduced to the cases $f\simeq^1 g$ and $f\simeq^{-1}g$. Assume
that $f\s^1g$ or $f\s^{-1}g$ by a shrinking map $h\colon J_{m}\rightarrow J_{n}$. Here we can always assume
that $m-n$ is an even integer. For if not, we can construct a shrinking map $\overline{h}\colon J_{m+1}\rightarrow J_{n}$ such
that $\overline{h}(m+1)=n$ and a path $\overline{f}\colon J_{m+1}\rightarrow G$ such that $\overline{f}\s^1 g$ or $\overline{f}\s^{-1} g$, and then obtain a digraph map $\overline{F}\colon C_{\overline{h}}\rightarrow G$ or $\overline{F}\colon C^{-}_{\overline{h}}\rightarrow G$.

\begin{description}
\item[\emph{Case 1}] If $f\simeq^{1} g$, then there is a digraph map $F\colon  C_{h}\rightarrow G$. Suppose the first slant line $(0,i)\rightarrow (1,k_0)$ of $C_h$ from left to right is in position $k_{0}$, that is, $h(i) = k_0$, $0\leq k_0< i$. Then we construct a subdivision $\widetilde{g}$ of $g$ by a shrinking map $H$, where $H\colon  J_{n+2}\rightarrow J_n$ is defined by
$$ H(k)=\left\{
                     \begin{array}{ll}
                       k, & \hbox{$k\leq k_{0}$;} \\
                       k-1, & \hbox{$k=k_{0}+1$;} \\
                       k-2, & \hbox{$k_{0}+2 \leq k$.}
                     \end{array}
                   \right
.$$

Next, we construct a new shrinking map $\widetilde{h}\colon  J_{m}\rightarrow J_{n+2}$ by $$\widetilde{h}(i) = \left\{
                     \begin{array}{ll}
                       i, & \hbox{$i\leq k_{0}+2$;} \\
                       h(i)+2, & \hbox{$k_{0}+3 \leq i \leq m $,}
                     \end{array}
                   \right
.$$  and a map $\widetilde{F}\colon C_{\widetilde{h}} \rightarrow G$ by
 $$\widetilde{F}((0, i)\rightarrow (1, j)) = \left\{
                                            \begin{array}{ll}
                                              F((0, i)\rightarrow (1, j)), & \hbox{$j\leq k_{0}$;} \\
                                              F((0, i)\rightarrow (1, j-1)), & \hbox{$j = k_{0}+1$;} \\
                                              F((0, i)\rightarrow (1, j-2)), & \hbox{$  k_{0}+2 \leq j\leq n+2$.}
                                            \end{array}
                                          \right.
 $$ Then $\widetilde{F}|_{0\Box J_{m}} = f$, $\widetilde{F}|_{1\Box J_{n+2}} = \widetilde{g}$, and $\widetilde{F}$ is a digraph map. Thus $f\simeq^{1}\widetilde{g}$.

\item[\emph{Case 2}]If $f\simeq^{-1} g$, there is a digraph map $F\colon  C^{-}_{h}\rightarrow G$. As in Case 1, we construct the same shrinking maps $\widetilde{h}$ and $H$. Further we define a map $\widetilde{F}\colon C^{-}_{\widetilde{h}} \rightarrow G$ by $$\widetilde{F}((1,j)\rightarrow (0,i)) = \left\{
                                            \begin{array}{ll}
                                               F((1,j)\rightarrow (0,i)), & \hbox{$j\leq k_{0}$;} \\
                                              F( (1,j-1)\rightarrow (0,i)), & \hbox{$j = k_{0}+1$;} \\
                                              F((1,j-2)\rightarrow (0,i)), & \hbox{$k_{0}+2\leq j\leq n+2 $.}
                                            \end{array}
                                          \right.
 $$
Then $\widetilde{F}|_{0\Box J_{m}} = f$, $\widetilde{F}|_{1\Box J_{n+2}} = \widetilde{g}$, and $\widetilde{F}$ is a digraph map. Thus $f\simeq^{-1}\widetilde{g}$.
\end{description}

In either case, as $m-n-2$ is an even number, we can iterate by considering $f\simeq^{1}\widetilde{g}$ or
$f\simeq^{-1}\widetilde{g}$ and the cylinder $C_{\widetilde{h}}$ or $C^{-}_{\widetilde{h}}$ of the
map $\widetilde{h}$. Repeat this process until $\widetilde{h}(m)=m$, in which case $\widetilde{h}$ and $\widetilde{g}$ are
denoted by $\overline{h}$ and $\overline{g}$ respectively. Then $f\dr\overline{g}$ or $f\dl\overline{g}$.
By construction, $\overline{g}$ is a subdivision of $g$. Hence $f\simeq_{1}g$ or $f\simeq_{-1}g$.
\end{proof}
 \end{lemma}
\begin{remark}\label{onec}
Lemma~\ref{equ} shows that $C$-homotopy and $F$-homotopy are equivalent. However, this equivalence does not hold at the level of one-step homotopy. Specifically, from the proof of Lemma~\ref{equ}, we see that if $f \simeq^{1} g$, then $f \simeq_1 g$; and if $f \simeq^{-1} g$, then $f \simeq_{-1} g$. The converses, however, are not generally true. We present an illustrative counterexample below.
\begin{example}Let 
\begin{center}
\begin{tikzpicture}
\node[below, black] at (7.4,0.27) {$\phi:$};
   %\node[above, black] at (6.5,0.2) {h};
    \node[below, black] at (8,0) {0};
     \node[below, black] at (9,0) {1};
     \node[below, black] at (10,0) {2};
     
\filldraw (8,0) circle (.04)
(9,0) circle (.04);
\filldraw (10,0) circle (.04);
   \draw[->, thick] (8.1,0) -- (8.9,0);
   \draw[->, thick] (9.9,0) -- (9.1,0);
\node[below, black] at (11,0.2) {$\longrightarrow$};
\node[below, black] at (12,0) {$\ast$};
\filldraw (12,0) circle (.04);
\end{tikzpicture}
\end{center}
and 
\begin{center}
\begin{tikzpicture}
\node[below, black] at (7.4,0.27) {$\psi:$};
   %\node[above, black] at (6.5,0.2) {h};
    \node[below, black] at (8,0) {0};
     \node[below, black] at (9,0) {1};
     \node[below, black] at (10,0) {2};
     
\filldraw (8,0) circle (.04)
(9,0) circle (.04);
\filldraw (10,0) circle (.04);
   \draw[<-, thick] (8.1,0) -- (8.9,0);
   \draw[<-, thick] (9.9,0) -- (9.1,0);
\node[below, black] at (11,0.2) {$\longrightarrow$};
\node[below, black] at (12,0) {$\ast.$};
\filldraw (12,0) circle (.04);
\end{tikzpicture}
\end{center}
Clearly, $\phi\simeq_1 \psi.$ However, $\psi $ is not one-step $C$-homotopic to $\psi$ since there is no shrinking map $h: J_2 \rightarrow I_2$. Suppose there is a shrinking map $h: J_2 \rightarrow I_2$, then $h(0)=0$ and $h(2) =2$, implying that $h(1) =1$ and $h$ is not a digraph map. This is a conflict. Thus $\psi$ is not one-step direct $C$-homotopic to $\psi$.
\end{example}
This example shows that, even though two trivial paths may be 
$C$-homotopic to each other, one may not be one-step 
$C$-homotopic to the other. This highlights the crucial role that arrow directions play in the homotopy theory of digraphs, as well as the inherent subtlety and mystery of the theory itself.
\end{remark}

In particular, if we focus on loops, $f\s^C g$ if and only if  $f\simeq_{F}g$. In what follows, when we refer to $C$-homotopy, we will use Lemma \ref{equ}. For simplicity of
notation, we will write $[J,G]$ for $[J,G]^C$. 
\section{ Homotopy Groups of Digraphs}
The most important goal of this section is to define the homotopy groups $\overline{\pi}_{n}(G)$ of a digraph $G$, which admit grid descriptions and are variations of the GLMY homotopy groups introduced in \cite{Y-homo}. And our homotopy groups of digraphs generalize of the $A$-homotopy groups of graphs in \cite{A-homo} and cubical sets in \cite{Carranza}. Furthermore, we prove some properties of homotopy
groups of digraphs that are similar to those of topological spaces. In particular,
we introduce a new reduced loop-digraph $\overline{L}G$ with the property that
$\overline{\pi}_n(\overline{L}G)\approx\overline{\pi}_{n+1}(G)$.

Recall the definition of the homotopy groups of a topological space $X$ with basepoint $x_0$. They can be defined by relative maps from the pair $(I^n,\partial I^n)$ to $(X,x_0)$, where $I^n$ is the $n$-dimensional cube and $\partial I^n$ is the
boundary of $I^n$. Following this idea, the $n$-dimensional homotopy group of a based
digraph $G$ with base-point $\ast$ is given by relative digraph maps from the $n$-dimensional grid digraph to $(G,\ast)$.
In this section we always assume that $n\geq1$ unless otherwise stated.
To picture the $n$-dimensional grid digraph vividly, it is necessary to introduce the relative box
product of digraph pairs.
\begin{definition}
Let $(G,A)$ and $(H,B)$ be digraph pairs. The \emph{relative box product} $(G,A)\Box(H,B)$ of $(G,A)$ and $(H,B)$ is the digraph pair $(G\Box H, A\Box H \cup G\Box B)$.
\end{definition}
\begin{example}By definition of the relative box product,
$$(J_3,\partial J_3)\Box(J_2,\partial J_2) = (J_3\Box J_2,\partial J_3\Box J_2 \cup J_3\Box \partial J_2 ),$$
which is shown as follows.
$$
\xymatrix@R=2.5em{
(0,2) \ar@{.>}[d] \ar@{.>}[r]  &(1,2)\ar[d] & (2,2)  \ar@{.>}[l] \ar@{.>}[r]  \ar[d]& (3,2)\ar@{.>}[d] \\
(0,1) \ar[r]        & (1,1)        & (2,1)  \ar[l]  \ar[r]       & (3,1)      \\
(0,0) \ar@{.>}[r]\ar@{.>}[u]  & (1,0)\ar[u]  & (2,0)  \ar@{.>}[l]\ar[u]\ar@{.>}[r]   & (3,0) \ar@{.>}[u]
, }
 $$
Here the dotted arrows form the sub-digraph $\partial J_3\Box J_2 \cup J_3\Box \partial J_2 $ of $J_3\Box J_2$. 
\end{example}
In general, the relative box product of two line digraphs
 is the $2$-dimensional grid digraph relative to its boundary digraph. More generally, if we consider the $n$-times relative box
product of line digraphs, it is exactly the $n$-dimensional grid digraph relative to its boundary digraph. With this in mind, we proceed to consider a new homotopy between digraph maps $f\colon  (I_{m_{i}},\partial I_{m_i})^{\Box n}\rightarrow (G,\ast)$ and $g\colon (I_{n_{i}},\partial I_{n_i})^{\Box n}\rightarrow (G,\ast)$ by using subdivision. Here, $(I_{m_{i}},\partial I_{m_i})^{\Box n}$ represents the relative box product $$(I_{m_{1}},\partial I_{m_1})\Box (I_{m_{2}},\partial I_{m_2})\Box \cdots \Box (I_{m_{n}},\partial I_{m_n}),$$ where each $m_{i}$ is the length of line digraph $I_{m_i}$, for all $1\leq i\leq n$. We will denote $(G,A)\Box (H,\emptyset)$ by $(G,A)\Box H.$

To define \( F \)-homotopy between \( n \)-dimensional grids, we first introduce the concepts of \( n \)-dimensional shrinking maps and subdivisions. A relative digraph map $
h\colon (I_{M_i}, \partial I_{M_i})^{\Box n} \to (I_{m_i}, \partial I_{m_i})^{\Box n}$
is called an \emph{\( n \)-dimensional shrinking map} if it is the box product of 1-dimensional shrinking maps, i.e.,
$
h = h_1 \Box h_2 \Box \cdots \Box h_n,
$
where each 
$
h_i\colon (I_{M_i}, \partial I_{M_i}) \to (I_{m_i}, \partial I_{m_i})$
is a shrinking map for \( 1 \leq i \leq n \). Given a relative digraph map 
$
f\colon (I_{m_i}, \partial I_{m_i})^{\Box n} \to (G, *),$
if there exists a shrinking map \( h \) as above, then the composition \( \overline{f} = f \circ h \) is called a \emph{subdivision} of \( f \). When \( n = 1 \), this definition coincides with the shrinking map and subdivision introduced in Section~2. For simplicity, we will omit the dimension and refer to all such maps as shrinking maps.

\begin{definition}\label{Fhomo}
Let $f\colon  (I_{m_{i}},\partial I_{m_i})^{\Box n}\rightarrow (G,\ast)$ and $g\colon (I_{n_{i}},\partial I_{n_i})^{\Box n}\rightarrow (G,\ast)$ be relative digraph maps. We call $f $ \emph{one-step $F$-homotopic} to $g$ if there exist subdivisions $\overline{f}$ and $\overline{g}$ of $f$ and $g$ respectively such that $f\dr g$, denoted by $f\simeq_{1}g$. We also call $g$ \emph{one-step inverse $F$-homotopic} to $f$ and write $g\s_{-1}f$.

 More generally, we say that $f$ is \emph{$F$-homotopic} to $g$ if there is a finite sequence of $\{f_i\}_{i=0}^l$ such that $f_{0}=f$, $f_{l}=g$ and there are one-step $F$-homotopies $f_i\s_{1}f_{i+1}$ or $f_i\s_{-1}f_{i+1}$ for $0\leq i\leq l-1$.
\end{definition}
Clearly, when $n=1$ Definition \ref{Fhomo} is same as the Definition \ref{1Fhomo}.  On the other hand, if there is some $i$ such that $m_i = 0$, then $(I_{m_{i}},\partial I_{m_i})^{\Box n}$ is an $(n-1)$-dimensional
 grid, not an $n$-dimensional grid. Suppose $m_1 =0$, one can easily verify that for any relative digraph map
$f\colon (I_{m_{i}},\partial I_{m_i})^{\Box n}\rightarrow (G,\ast) $, there is a subdivision
$\overline{f}\colon (I_{1},\partial I_{1})\Box(I_{m_i},$ $\partial I_{m_i})^{\Box (n-1)}\rightarrow (G,\ast)$ of
$f$ such that $\overline{f}\s_1 f.$ As the homotopy type is not affected, we always assume that $m_i\geq 1$
for $1\leq i\leq n$.

Obviously, $F$-homotopy is an equivalence relation on 
$$\bigsqcup\limits_{I_{m_{i}},\forall i}\operatorname{Hom}((I_{m_{i}},\partial I_{m_i})^{\Box n};(G,\ast)).$$ 
Fix a digraph map $f\colon (I_{m_{i}},\partial I_{m_{i}})^{\Box (n+1)}\rightarrow (G,\ast)$. Any subdivision of $f$ can be decomposed into a series of coordinate subdivisions, just like coordinate projections. Assume that $\overline{f}$ is the subdivision of $f$ given by $q = h_1\Box h_2\Box \cdots \Box h_{n+1}$. Then $q$ can be decomposed as the composition of $h_1\Box h_2\Box \cdots \Box h_{n}\Box id $ and $id\Box id\Box \cdots \Box id\Box h_{n+1} $ shown as follows:
$$\begin{small}
\xymatrix{
  (I_{M_i},\partial I_{M_i})^{\Box n} \Box (I_{m_{n+1}},\partial I_{m_{n+1}})  \ar[rr]^{h_1\Box h_2\Box \cdots \Box h_n \Box id } & &  (I_{m_i},\partial I_{m_i})^{\Box n} \Box (I_{m_{n+1}},\partial I_{m_{n+1}}) \ar[r]^{\quad \quad \quad \quad \quad \quad  f} &  (G,\ast)     \\
 (I_{M_i},\partial I_{M_i})^{\Box n} \Box (I_{M_{n+1}},\partial I_{M_{n+1}}) \ar[u]^{id\Box id\Box \cdots\Box id \Box h_{n+1} }. \ar[urr]^{q} \ar[urrr]_{\overline{f}} }
 \end{small}
$$
Moreover we can further decompose $q$ one coordinate at a time, and such a decomposition is independent of the order of coordinate subdivisions.

 \begin{remark}\label{nstan} By Lemma \ref{standard} and the above decomposition, for any digraph map $f\colon (I_{m_i},\partial I_{m_i})^{\Box n} \rightarrow (G,\ast)$, there exists a subdivision $\overline{f}\colon (J_{M_i},\partial J_{M_i})^{\Box n} \rightarrow (G,\ast)$ of $f$ such that $\overline{f}\s_1 f.$ Hereafter, we always consider $f$ as a digraph map from the standard $n$-grid and call this $f$ the standard $n$-grid map.

 It is well-known that $\mathbb{N}^{\times n}$ determines a poset $(\mathbb{N}^{\times n},\leq)$, where  $(m_1,m_2,\cdots,m_n)\leq (s_1,s_2,\cdots,s_n)$ if and only if  $m_i\leq s_i$ for all $i$. Similar to the directed system $\{\operatorname{Hom}((J_{m},$ $\partial J_{m});(G,\ast)); l_{m}^{n}\}_{m\leq n}$, there is a directed system $\{\operatorname{Hom}((J_{m_i},\partial J_{m_i})^{\Box n};(G,\ast)); l_{M}^{S}\}_{M\leq S}$ based on the directed set $(\mathbb{N}^{\times n},\leq)$, where $$l_{M}^{S}\colon \operatorname{Hom}((J_{m_i},\partial J_{m_i})^{\Box n};(G,\ast))\rightarrow \operatorname{Hom}((J_{s_i},\partial J_{s_i})^{\Box n};(G,\ast)),\text{ }f\mapsto \overline{f},$$ with $\overline{f}(i_1,i_2,\cdots,i_n) =\left\{
                     \begin{array}{ll}
                       f(i_1,i_2,\cdots,i_n), & \hbox{$(i_1,i_2,\cdots,i_n)\leq M$;} \\
                       f(m_1,m_2,\dots,m_n), & \hbox{otherwise.}
                     \end{array}
                   \right
.$

The direct limit of $\{\operatorname{Hom}((J_{m_i},\partial J_{m_i})^{\Box n};(G,\ast)); l_{M}^{S}\}_{M\leq S}$ is
$$\lim\limits_{\rightarrow}\operatorname{Hom}((J_{m_i},\partial J_{m_i})^{\Box n};(G,\ast))= \bigsqcup\limits_{m_i \geq1, \forall i}\operatorname{Hom}((J_{m_i},\partial J_{m_i})^{\Box n};(G,\ast))/\sim,$$ where $f_{M}\sim f_{S}$ if and only if there exists $V\in \mathbb{N}^{\times n}$ with $M\leq V$ and $S\leq V$ such that $l_{M}^{V}(f_{M}) = l_{S}^{V}(f_{S})$ for any $f_{M}\in \operatorname{Hom}((J_{m_i},\partial J_{m_i})^{\Box n};(G,\ast))$ and $f_{S}\in \operatorname{Hom}((J_{s_i},$ $\partial J_{s_i})^{\Box n};(G,\ast))$. From the definition of direct limit, there is a family of injections $$\{i_M\colon \operatorname{Hom}((J_{m_i},\partial J_{m_i})^{\Box n};(G,\ast))\rightarrow \lim\limits_{\rightarrow}\operatorname{Hom}((J_{m_i},\partial J_{m_i})^{\Box n};(G,\ast)) \}$$ such that $i_S\circ l_{M}^{S} = i_M.$

For the quotient map $$\Gamma\colon \bigsqcup\limits_{M=(m_1,\cdots,m_n)}\operatorname{Hom}((J_{m_i},\partial J_{m_i})^{\Box n};(G,\ast))\rightarrow \lim\limits_{\rightarrow}\operatorname{Hom}((J_{m_i},\partial J_{m_i})^{\Box n};(G,\ast))$$
and elements $\{f\},\text{ }\{g\}\in \lim\limits_{\rightarrow}\operatorname{Hom}((J_{m_i},\partial J_{m_i})^{\Box n};(G,\ast))$,
we write $\{f\}\s_{1}\{g\}$ if and only if
$f\s_1 g$, $i_M(f) = \{f\}$ and $i_{S}(g)=\{g\}$. One can easily check that this is independent of the representatives $M$ and $S$. To simplify notation, we shall denote $\lim\limits_{\rightarrow}\operatorname{Hom}((J_{m_i},\partial J_{m_i})^{\Box n};(G,\ast))$ by $\operatorname{Hom}((J,\partial J)^{\Box n};(G,\ast))$ and the set of $F$-homotopy classes of $\operatorname{Hom}((J,\partial J)^{\Box n};(G,\ast))$ by $[(J,\partial J)^{\Box n};(G,\ast)]$.

In fact, $[(J,\partial J)^{\Box n};(G,\ast)]$ is the underlying set of what will be our $n$-dimensional homotopy group. Before we give the group structure of homotopy group of digraphs, the following important property of subdivision should be considered. This is because we need to take into account the arrows direction in the digraphs.
\begin{lemma}\label{subre}
Let $f\colon (J_{m},\partial J_{m})\rightarrow (G,\ast)$ be a standard loop. For any two subdivisions
$$f_1\colon (J_{m_1},\partial J_{m_1})\rightarrow (G,\ast)$$ and $$f_2\colon (J_{m_2},\partial J_{m_2})\rightarrow (G,\ast)$$
of $f$, there exists a common subdivision $\overline{f}$ of $f_1$ and $f_2$.
\begin{proof}
 Suppose that there are subdivisions 
 $$h_1\colon (J_{m_1},\partial J_{m_1})\rightarrow (J_{m},\partial J_{m})~~ \text{and}~~h_2\colon (J_{m_2},\partial J_{m_2})\rightarrow (J_{m},\partial J_{m}) $$ such that $f_1 = f\circ h_1$ and $f_2 = f\circ h_2$.
To prove the lemma, it is sufficient to construct two shrinking maps $q_1\colon (J_{M},\partial J_{M})\rightarrow (J_{m_1},\partial J_{m_1}) $ and $q_2\colon (J_{M},\partial J_{M})\rightarrow (J_{m_2},\partial J_{m_2}) $ such that $h_1\circ q_1 = h_2\circ q_2$, giving a commutative diagram
$$
\xymatrix{
  (J_{M},\partial J_{M}) \ar[d]_{q_2} \ar[r]^{q_1}
                & (J_{m_1},\partial J_{m_1}) \ar[d]^{h_1}  \\
  (J_{m_2},\partial J_{m_2}) \ar[r]_{h_2}
                & (J_{m},\partial J_{m})\ar[r]^{f} & (G,\ast).             }
$$
To start, let $i=j_1=j_2 = j =0$ and $q_0(0)=q_1(0)=0$.

Fix $i$. Denote $N = max\{\mid h_1^{-1}(i)\mid-1,\mid h_2^{-1}(i)\mid-1\}$ and $l = 1.$
\begin{enumerate}
\item If $\mid h_1^{-1}(i)\mid> \mid h_2^{-1}(i)\mid$, let $q_1(j+l) =  j_1+l $ and if $l<  \mid h_2^{-1}(i)\mid$, let $q_2(j+l) = j_2+l$, otherwise $q_2(j+l) = j_2+\mid h_2^{-1}(i)\mid-1;$
\item If $\mid h_1^{-1}(i)\mid = \mid h_2^{-1}(i)\mid$, let $q_1(j+l) =  j_1+l $ and $q_2(j+l) = j_2+l$;
\item If $\mid h_1^{-1}(i)\mid< \mid h_2^{-1}(i)\mid$, let $q_2(j+l) =  j_2+l $ and if $l<  \mid h_1^{-1}(i)\mid$, let $q_1(j+l) = j_1+l$, otherwise $q_1(j+l) = j_1+\mid h_1^{-1}(i)\mid-1.$
\end{enumerate}
If $l=N +1$, we stop and let $j=j+N$, $j_1 = j_1+\mid h_1^{-1}(i)\mid-1$, $j_2 = j_2+\mid h_2^{-1}(i)\mid-1$. If $l\leq N$, increase $l$ by 1. Then we iterate the above process. If $i\leq m-1$, increase $i$ by 1 and iterate the above process, otherwise we stop.

Then let $M = j$, $q_1\colon  J_M\rightarrow J_{m_1}$ and  $q_2\colon  J_M\rightarrow J_{m_2}$. This gives what we need to construct.

Now let us check $q_1$ is a digraph map. For any arrow of $j\rightarrow j+1$, we need to verify $q_1(j\rightarrow j+1)$ is a vertex or an arrow $q_1(j)\rightarrow q_1(j+1)$. By the above process, there is always an arrow between $q_1(j)$ and $q_1(j+1)$ or $q_1(j)= q_1(j+1)$. We only need to check if $q_1(j)\neq q_1(j+1)$, then $q_1(j)\rightarrow q_1(j+1)$.

There is a very interesting observation. For any $i< m$, $\left| \left| h_1^{-1}(i)\right|- \left| h_2^{-1}(i)\right| \right|$ must be even. If not, there exists an integer $i<m$ such that $\left| \left| h_1^{-1}(i)\right|- \left| h_2^{-1}(i)\right| \right|$ is odd, implying that the arrow direction connecting the last vertex $k_1$ of $h_1^{-1}(i)$ with the first vertex $k_1+1$ of $h_1^{-1}(i+1)$ in $J_{m_1}$ is different from  the arrow direction connecting the last vertex $k_2$ of $h_2^{-1}(i)$ with the first vertex $k_2+1$ of $h_2^{-1}(i+1)$ in $J_{m_2}$. Assume that $i\rightarrow i+1$ in $J_m$, then  at least one of $h_1$ and $h_2$ is not a digraph map. This is a contradiction. Hence for any $i< m$, $\left| \left| h_1^{-1}(i)\right|- \left| h_2^{-1}(i)\right| \right|$ must be even. So by the process of constructing $q_1$, $\left| max-\left| h_1^{-1}(i)\right|+1 \right|-\left| max-\left| h_2^{-1}(i)\right|+1 \right|$ must be even, therefore $q_1(j)\rightarrow q_1(j+1)$. Thus $q_1|_{q_1^{-1}([0,m-1])}$ is a digraph map. Finally, we consider the last vertex $i=m$. It is clear that $q_1$ is a digraph map by the iterated construction for $l$.  Similar to $q_1$, the map $q_2$ is also a digraph map.

It is clear that $q_1$ and $q_2$ are surjective, preserve order, and satisfy $h_1\circ q_1 = h_2\circ q_2$, that is,
$q_1$ and $q_2$ are shrinking maps. Therefore $\overline{f}=f \circ h_1\circ q_1$ is exactly the common subdivision of $f_1$ and $f_2$.
\end{proof}
\end{lemma}
\end{remark}

Furthermore, for a higher dimensional standard grid map $f\colon (J_{m_i},\partial J_{m_i})^{\Box n}\rightarrow (G,\ast)$, we immediatly have the same result.
\begin{corollary}\label{subren}
Let $f\colon (J_{m_i},\partial J_{m_i})^{\Box n}\rightarrow (G,\ast)$. For any two subdivisions
$f_1$ and
$f_2$ of $f$, there is a common subdivision $$\overline{f}\colon (J_{N_i},\partial J_{N_i})^{\Box n}\rightarrow (G,\ast)$$ of $f_1$ and $f_2$. \qed

\end{corollary}
Now let us give the group structure of $[(J,\partial J)^{\Box n};(G,\ast)].$
\begin{proposition}\label{group}Let $G$ be a based digraph. Then $[(J,\partial J)^{\Box n};(G,\ast)]$ is a group for $n\geq1$.
\begin{proof}
To give the group structure of $[(J,\partial J)^{\Box n};(G,\ast)]$, we divide the proof into three steps.
\begin{description}
    \item[\emph{Step 1}] \emph{First, let us define the multiplication more generally in 
    $$\bigsqcup\limits_{I_{m_{i}}\atop 0\leq i \leq n}\operatorname{Hom}((I_{m_{i}}, \partial I_{m_{i}})^{\Box n};(G,\ast)).$$} In fact, there are several multiplications along different coordinates, so we can define a family of multiplications $\{\mu^j\}_{j=1}^{n}$, where 
    $$
 \bigsqcup\limits_{I_{m_{i}} \atop 0\leq i \leq n}\operatorname{Hom}((I_{m_{i}},\partial I_{m_{i}})^{\Box n};(G,\ast))\times\bigsqcup\limits_{I_{m_{i}}\atop 0\leq i \leq n}\operatorname{Hom}((I_{m_{i}},\partial I_{m_{i}})^{\Box n};(G,\ast))$$
$$ \downarrow \mu^j$$
$$\bigsqcup\limits_{I_{m_{i}}\atop 0\leq i \leq n}\operatorname{Hom}((I_{m_{i}},\partial I_{m_{i}})^{\Box n};(G,\ast)).
$$
    
    % $$\tiny{\bigsqcup\limits_{I_{m_{i}} \atop 0\leq i \leq n}\operatorname{Hom}((I_{m_{i}},\partial I_{m_{i}})^{\Box n};(G,\ast))\times\bigsqcup\limits_{I_{m_{i}}\atop 0\leq i \leq n}\operatorname{Hom}((I_{m_{i}},\partial I_{m_{i}})^{\Box n};(G,\ast))\stackrel{\mu^j}{\rightarrow} \bigsqcup\limits_{I_{m_{i}}\atop 0\leq i \leq n}\operatorname{Hom}((I_{m_{i}},\partial I_{m_{i}})^{\Box n};(G,\ast)).}$$
Suppose $f\colon (I_{m_{i}},\partial I_{m_{i}})^{\Box n}\rightarrow (G,\ast)$ and $g\colon  (I_{n_{i}},\partial I_{n_{i}})^{\Box n}\rightarrow (G,\ast)$. Let $M_i = max\{n_i,$ $ m_i\}$ for $i\ne j$. 
To ensure that the domain of \( \mu^j(f, g) \) remains a grid, we extend the domains of \( f \) and \( g \) to the grids
\[
\widetilde{f} \colon (I_{M_1}, \partial I_{M_1}) \Box \cdots \Box (I_{m_j}, \partial I_{m_j}) \Box \cdots \Box (I_{M_n}, \partial I_{M_n}) \to (G, \ast)
\]
and
\[
\widetilde{g} \colon (I_{M_1}, \partial I_{M_1}) \Box \cdots \Box (I_{n_j}, \partial I_{n_j}) \Box \cdots \Box (I_{M_n}, \partial I_{M_n}) \to (G, \ast)
\]
respectively. The idea is straightforward shown as follows.
$$
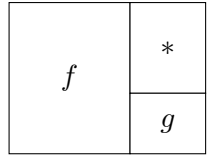
\begin{tikzpicture}
\node (rect) at (0,0) [draw, thin, fill=white, fill opacity=0.7,
minimum width=1.6cm,minimum height=2cm] {$f$};
\node (rect) at (1.3,0.4) [draw, thin, fill=white, fill opacity=0.7,
minimum width=1cm,minimum height=1.2cm] {$\ast$};
\node (rect) at (1.3,-0.6) [draw, thin, fill=white, fill opacity=0.7,
minimum width=1cm,minimum height=0.8cm] {$g$};
\end{tikzpicture}$$

Then we define $$ \mu^j(f,g)\colon  (I_{M_{1}},\partial I_{M_{1}})\Box\cdots \Box (I_{(m_{j}+n_j)},\partial I_{(m_{j}+n_j)})\Box\cdots \Box(I_{M_{n}},\partial I_{M_{n}}) \rightarrow (G,\ast)$$
by
\begin{align*}  \mu^j(f,g)(i_1,i_2,...,i_j,...,i_n)&=(\widetilde{f} \vee \widetilde{g})(i_1,i_2,...,i_j,...,i_n) \\ &=\left\{
                     \begin{array}{ll}
                       \widetilde{f}(i_1,i_2,...,i_j,...,i_n), & \hbox{$i_j\leq m_j$;} \\
                        \widetilde{g}(i_1,i_2,...,i_j-m_j,...,i_n), & \hbox{$i_j> m_j$.}
                     \end{array}
                   \right
                   .
\end{align*}
Next let us check that $\mu^j(f,g)\s_1 \mu^j(f^{'},g)$ if $f\s_1 f^{'}$. Suppose
$\mu^j(f,g)=\widetilde{f} \vee \widetilde{g} $, $\mu^j(f',g)=\widetilde{f'} \vee \widetilde{g}'$, and there exist subdivisions $\overline{f}$ and $ \overline{f^{'}}$ of $f$ and $f^{'}$ respectively
such that $\overline{f}\rightrightarrows \overline{f^{'}}$ by $h$ and $h^{'}$. By Corollary \ref{subren},
we have a common subdivision $\overline{\widetilde{f}}$ of $\widetilde{f}$ and $\overline{f}$ and a
common subdivision $ \overline{\widetilde{f^{'}}}$ of $\widetilde{f^{'}}$ and $\overline{f^{'}}$ such that $\overline{\widetilde{f}}\dr \overline{\widetilde{f^{'}}}$. Then clearly we have $\mu^j(\overline{\widetilde{f}},g)\dr \mu^j(\overline{\widetilde{f^{'}}},g)$, therefore  $\mu^j(\overline{\widetilde{f}},g)$ and $\mu^j(\overline{\widetilde{f'}},g)$ are subdivisions of $\mu^j(f,g)$ and $\mu^j(f',g)$ respectively, so $\mu^j(f,g)\s_1 \mu^j(f^{'},g)$. Similarly, if $g\s_1 g^{'}$, then $\mu^j(f,g)\s_1 \mu^j(f,g^{'})$. Moreover, $\mu^j(f,g)\s_1 \mu^j(f,g^{'})\s_1 \mu^j(f^{'},g^{'})$. Hence  if $f\s_F f^{'}$ or $g\s_F g^{'}$, then $\mu^j(f,g)\s_F \mu^j(f^{'},g^{'})$.
Subsequently, we default to using this standard $n$-grid map $f$ with even length as a representation for each $\{f\}$ in $\lim\limits_{\rightarrow}\operatorname{Hom}((J_{m_i},\partial J_{m_i})^{\Box n};(G,\ast))$. Further,
to prove $\{f\}\s_{1} \{g\}$ it suffices to prove $f\s_1 g$ for any elements $\{f\},\text{ }\{g\}\in \lim\limits_{\rightarrow}\operatorname{Hom}((J_{m_i},\partial J_{m_i})^{\Box n};(G,\ast))$. Then the digraph map $\mu^j$ induces a map $$\mu^j_{n}\colon [(J,\partial J)^{\Box n};(G,\ast)]\times[(J,\partial J)^{\Box n};(G,\ast)]\rightarrow [(J,\partial J)^{\Box n};(G,\ast)]$$ with  $\mu^j_{n}([f],[g])\mapsto [\mu^j(f,g)]$. By the previous paragraph, $\mu^j_{n}$ is well-defined.

 \item[\emph{Step 2}] \emph{$\mu^j_{n}([f],[g]) = \mu^i_{n}([f],[g]).$}  The proof is straightforward, and the idea is illustrated in the following diagram.
%By Step 1, $\widetilde{f}$ has the same $F$-homotopy type as $f$, and $\widetilde{g}$ has the same $F$-homotopy type as $g$. Then we have the following picture under the $F$-homotopy equivalence relation.

\[
\mu_n^i(f, g) =
\tikz[baseline=(current bounding box.center)]{
\node (rect) at (0,0) [draw, thin, fill=white, fill opacity=0.7,
minimum width=0.6cm,minimum height=1cm] {$f$};
\node (rect) at (0.8,0.3) [draw, thin, fill=white, fill opacity=0.7,
minimum width=1cm,minimum height=0.4cm] {$\ast$};
\node (rect) at (0.8,-0.2) [draw, thin, fill=white, fill opacity=0.7,
minimum width=1cm,minimum height=0.6cm] {$g$};
}
=
\tikz[baseline=(current bounding box.center)]{
 \node (rect) at (0,0) [draw, thin, fill=white, fill opacity=0.7,
minimum width=0.6cm,minimum height=1cm] {$f$};
\node (rect) at (0.8,0.3) [draw, thin, fill=white, fill opacity=0.7,
minimum width=1cm,minimum height=0.4cm] {$\ast$};
\node (rect) at (0.8,-0.2) [draw, thin, fill=white, fill opacity=0.7,
minimum width=1cm,minimum height=0.6cm] {$g$};
\node (rect) at (0,1) [draw, thin, fill=white, fill opacity=0.7,
minimum width=0.6cm,minimum height=1cm] {$\ast$};
\node (rect) at (0.8,1) [draw, thin, fill=white, fill opacity=0.7,
minimum width=1cm,minimum height=1cm] {$\ast$};
}
=
\tikz[baseline=(current bounding box.center)]{
  \node (rect) at (0,0) [draw, thin, fill=white, fill opacity=0.7,
minimum width=0.6cm,minimum height=1cm] {$f$};
\node (rect) at (0.8,1.3) [draw, thin, fill=white, fill opacity=0.7,
minimum width=1cm,minimum height=0.4cm] {$\ast$};
\node (rect) at (0,1) [draw, thin, fill=white, fill opacity=0.7,
minimum width=0.6cm,minimum height=1cm] {$\ast$};
\node (rect) at (0.8,0) [draw, thin, fill=white, fill opacity=0.7,
minimum width=1cm,minimum height=1cm] {$\ast$};
\node (rect) at (0.8,0.8) [draw, thin, fill=white, fill opacity=0.7,
minimum width=1cm,minimum height=0.6cm] {$g$};
}
=
\tikz[baseline=(current bounding box.center)]{
 \node (rect) at (0,0) [draw, thin, fill=white, fill opacity=0.7,
minimum width=0.6cm,minimum height=1cm] {$f$};
\node (rect) at (0.2,0.8) [draw, thin, fill=white, fill opacity=0.7,
minimum width=1cm,minimum height=0.6cm] {$g$};
\node (rect) at (0.475,0) [draw, thin, fill=white, fill opacity=0.7,
minimum width=0.2cm,minimum height=1cm] {$\ast$};
}
= \mu_n^j(f, g)
\]
\item[\emph{Step 3}] \emph{$[(J,\partial J)^{\Box n};(G,\ast)]$ with the multiplication $\mu^j_{n}$ forms a group.}
It is easily seen that the multiplication is associative and the constant loop $e\colon (J_2,\partial J_2)^{\Box n}\rightarrow \ast$ is the unit loop.

Following the idea of the fundamental group in \cite[Lemma 4.19]{Y-homo}, the inverse $f^{-1}_j$ of $f$ along the $j$-th coordinate is defined by  $$f_j^{-1}((i_1,i_2,...,i_j,...,i_n))\mapsto f((i_1,i_2,...,l_j-i_j,...,i_n)).$$ And one easily prove that the inverse of $[f]$ is  independent of the coordinate $j$. Subsequently, we denote the inverse of $[f]$ by $[f^{-1}]$.
\end{description}
\end{proof}
\end{proposition}
Hereafter, the multiplication of $[(J,\partial J)^{\Box n};(G,\ast)]$ will be defined along the first coordinate. For simplicity, we will denote $\mu^1_{n}([f],[g])$ by $[f]\cdot[g]$ and $\mu^1(f,g)$ by $f\cdot g$.  With this group structure, we define the homotopy groups of digraphs.
\begin{definition}\label{ndef}
Let $G$ be a based digraph. The $n$-dimensional homotopy group $\overline{\pi}_n(G)$ is defined by $$\overline{\pi}_n(G):= [(J,\partial J)^{\Box n};(G,\ast)],\quad \quad \quad n\geq1.$$
\end{definition}
By \cite[Proposition 6.5]{Y-homo}, it is easy to prove that for any subdivision $\overline{f}$ of digraph map $f\colon  (J_{m_i},\partial J_{m_i})^{\Box n}\rightarrow (G,\ast)$ in double digraph category, $\mathcal{O}^{-1}(\overline{f})$ is $A$-homotopic to $\mathcal{O}^{-1}(f)$. Then Definition \ref{ndef} is the generalization of the $n$-dimensianal $A$-homotopy group $A_n^{G}(G)$ of graph $G$ in \cite{A-homo} with defining the $n$-grid map $f\colon  \mathbb{N}^n \rightarrow G$ and cubical set in \cite{Carranza}.

And by Theorem \ref{funda} and Lemma \ref{equ}, $\overline{\pi}_1(G)=\pi_{1}(G)$. It should be pointed out that whether our $n$-dimensional homotopy group is isomorphic to the $n$-dimensional homotopy group $\pi_n(G)$ introduced by Yau et al. is still unclear (see more details in Proposition \ref{0i}). So we use different notation $\overline{\pi}_n(G)$ instead of $\pi_n(G).$ For $n=0$, we use the definition of $\pi_0(G)$ introduced by Yau et al. \cite{Y-homo}, but for convenience we will denote this by $\overline{\pi}_0(G)$. Later we will give a more precise description about the relation between them. Analogously to the properties of homotopy groups of spaces, we obtain the following similar results.

\begin{proposition}\label{induce}
Each based digraph map $f\colon  (G,g_0)\rightarrow (H,h_0)$ induces a morphism $f_n\colon \overline{\pi}_n(G)\rightarrow \overline{\pi}_n(H)$ for all $n\geq0.$ If $n\geq 1$, $f_n$ is a homomorphism.
\begin{proof}
This is straightforward to prove and will be left to the reader.
\end{proof}
\end{proposition}
Clearly $\overline{\pi}_n$ is a functor from the category of based digraphs to the category of groups since $\overline{\pi}_n(\phi \circ \psi) = \overline{\pi}_n(\phi)\circ \overline{\pi}_n(\psi)$ and $\overline{\pi}_n (id_G) = id_{\overline{\pi}_n(G)}$ for all $n\geq 1$.
\begin{proposition}\label{map}
If $f\simeq g\colon  (G,g_0)\rightarrow (H,h_0)$, then $\overline{\pi}_n(f)= \overline{\pi}_n(g)$ for $n\geq0$.
\begin{proof}
If $n=0$, by \cite[Proposition 
 4.2]{Y-homo}, $\overline{\pi}_n(f)= \overline{\pi}_n(g)$.
Then we consider the case $n\geq 1$. If $f\simeq g\colon  (G,g_0)\rightarrow (H,h_0)$, then there is a digraph map $$F\colon  G\Box I_m\rightarrow H$$ such that $F|_{G\Box 0} =f$, $F|_{G\Box m} =g$ and $F|_{g_0\Box I_m} =h_0$. Assume $m=1$, that is, $f\dr g.$

For any homotopy class $[\gamma]$ of $\overline{\pi}_n(G)$, where $\gamma\colon  (J_{m_i},\partial J_{m_i})^{\Box n}\rightarrow (G,g_0)$, $ \overline{\pi}_n(f)([\gamma])$ $= [f\circ \gamma]$ and $\overline{\pi}_n(g)([\gamma])= [g\circ \gamma]$. Therefore we can construct a digraph map $F\circ (\gamma\Box id_{J_1})=\widetilde{F}\colon  (J_{m_i},\partial J_{m_i})^{\Box n}\Box(J_1,\emptyset)\rightarrow (H,h_0)$ such that $\widetilde{F}|_{(J_{m_i},\partial J_{m_i})^{\Box n}\Box 0}$ $=f\circ \gamma$ and $\widetilde{F}|_{(J_{m_i},\partial J_{m_i})^{\Box n}\Box1}=g\circ \gamma$. Thus $f\circ \gamma\dr g\circ \gamma$. More generally, $f\circ \gamma\s_F g\circ \gamma$ if $f\s g$.

\end{proof}
\end{proposition}
Following this proposition, one can easily check that the homotopy groups of digraphs are homotopy invariant. More precisely, we have the following corollary.
\begin{corollary}If $(G, g_0)\s (H , h_0)$, then $\overline{\pi}_n(G)\approx \overline{\pi}_n(H)$ for $n\geq 0.$ \qed
\end{corollary}
It is well-known that the homotopy group $\pi_n(X)$ for a based space $X$ is abelian if $n\geq2$. In our digraph version, we have the same conclusion. It is straightforward to prove.
\begin{proposition}\label{gp}
 If $n\geq 2$, then $\overline{\pi}_n(G)$ is abelian.  \qed
\end{proposition}
Due to the fact that $\pi_n(X)\approx \pi_{n-1}(\Omega X)$ for any space $X$, we try to obtain a similar property for digraphs. Based on the loop-digraph $LG$ defined in \cite{Y-homo} and the idea of simplicial homotopy theory in \cite{Curtis}, we define a reduced loop-digraph $\overline{L}G$ for based digraph $G$ with base-point $\ast$. Further, we prove that  $\overline{\pi}_{n}(G)\approx\overline{\pi}_{n-1}(\overline{L}G)$ for $n\geq1$.

Before defining the reduced loop-digraph, we first introduce the notion of a \emph{minimal path} and an equivalence relation on all paths in a based digraph \( G \), called \emph{subdivision equivalence}. The minimal path \( f_{\min} \) of a path \( f \colon I_m \to G \) is obtained by collapsing each vertex \( j+1 \) to \( j \) whenever \( f(j) = f(j+1) \). From the process of constructing of the minimal path, it can be seen as the inverse operation of subdivision. 
 It is easy to verify that \( f \) is one-step direct \( C \)-homotopic to \( f_{\min} \), that is, \( f \simeq^1 f_{\min} \). If two paths \( f \) and \( f' \) have the same minimal path, then we say that \( f \) and \( f' \) belong to the same subdivision class. Clearly, this defines an equivalence relation. We denote the subdivision class of a path \( f \) by \( \langle f \rangle \).

\begin{definition}
Let $G$ be a digraph with base-point $\ast$. The \emph{reduced path-digraph} $\overline{P}G$ is a based digraph  with base-point $\langle l_{\ast}\rangle$ whose vertex set consists of all subdivision classes $\langle f\rangle$ of paths on $G$, where $l_{\ast}\colon  (J_{2},\partial J_2)\rightarrow \ast$ and whose arrow set is defined by saying there is an arrow $\langle f\rangle\rightarrow \langle g\rangle$ provided there exist $f_1\in \langle f\rangle$ and $g_1\in\langle g\rangle$ such that $f_1\dr g_1$.
\end{definition}
\begin{definition}\label{rlp}
The \emph{reduced loop-digraph} $\overline{L}G$ is the based sub-digraph of $\overline{P}G$ whose vertex set consists of all subdivision classes $\langle f\rangle$ of loops on $G$.
\end{definition}

As mentioned, any path \( f \) can be subdivided into a standard path, and any representatives \( f \) and \( f' \) in the class \( \langle f \rangle \) are the subdivisions of the minimal path \( f_{\min} \). Combined with Lemma~\ref{subre}, we have for any \( f \in \langle f \rangle \), \( g \in \langle g \rangle \), and \( \langle f \rangle \rightarrow \langle g \rangle \), if \( f' \in \langle f \rangle \) and \( g' \in \langle g \rangle \), then there exist subdivisions \( \overline{\overline{\overline{f}}} \) of $f'$ and \(\overline{\overline{\overline{g}}} \) of $g'$ such that \( \overline{\overline{\overline{f}}} \dr\overline{\overline{\overline{g}}} \).
 Clearly this can be generalized to the  subdivision classes of $n$-dimensional gird maps $f\colon (J_{m_{i}},\partial J_{m_{i}})^{\Box n}\rightarrow (G,\ast)$ and $g\colon (J_{n_{i}},\partial J_{n_{i}})^{\Box n}\rightarrow (G,\ast)$, $n\geq 2.$ The proof can be illustrated by the following diagram:
$$
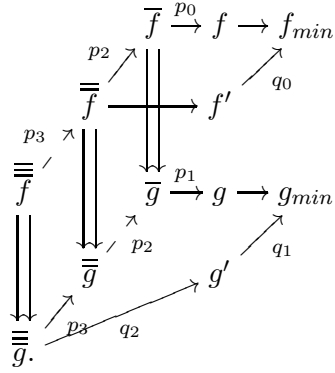
\xymatrix@R=2ex@C=2ex{
        &        & \overline{f}\ar@<.5ex>[dd] \ar@<-0.5ex>[dd]\ar[r]^{p_0} & f \ar[r] & f_{min} \\
& \overline{\overline{f}} \ar@<.5ex>[dd] \ar@<-0.5ex>[dd]\ar[rr]\ar[ur]^{p_2}
      & & f' \ar[ur]_{q_0}     &   \\
    \overline{\overline{\overline{f}}}  \ar@{-->}[ur]^{p_3} \ar@<.5ex>[dd] \ar@<-0.5ex>[dd] &          & \overline{g} \ar[r]^{p_1} & g \ar[r]&g_{min} \\
 & \overline{\overline{g}}\ar@{-->}[ur]_{p_2}
                &  & g'  \ar[ur]_{q_1} &    \\
\overline{\overline{\overline{g}}}.\ar[ur]_{p_3}\ar[urrr]_{q_2} & & & &}
$$
In fact, one can easily check that the homotopy type of $\overline{P}G$ is independent of the choice of base-point within the same path-component of $G$.
\begin{proposition}Let $\mathcal{DG^{\ast}}$ be the based digraph category whose objects are based digraphs and whose morphisms are digraph maps preserving the base-point. Then 
  $$\overline{L}\colon \mathcal{DG^{\ast}} \rightarrow \mathcal{DG^{\ast}},\quad G\to \overline{L}G $$ is a functor. 
\begin{proof}Let \( f \in \mathrm{Mor}(\mathcal{DG}^*) \). Clearly, for any loop \( \gamma \) in \( G \), the composition \( f \circ \gamma \) is a loop in \( H \). If \( \gamma \) and \( \gamma' \) are two representatives in the subdivision class \( \langle \gamma \rangle \), then both are subdivisions of the same minimal loop \( \gamma_{\min} \), via shrinking maps \( q \) and \( q' \). It follows that \( f \circ \gamma \) and \( f \circ \gamma' \) are subdivisions of the same loop in \( H \), so they lie in the same class: \( \langle f \circ \gamma \rangle = \langle f \circ \gamma' \rangle \). Therefore, the map
\[
\overline{L}f \colon \overline{L}G \to \overline{L}H, \quad \langle \gamma \rangle \mapsto \langle f \circ \gamma \rangle
\]
is well-defined. Moreover, if \( \langle \gamma \rangle \rightarrow \langle \eta \rangle \), then either \( \langle f \circ \gamma \rangle = \langle f \circ \eta \rangle \) or \( \langle f \circ \gamma \rangle \rightarrow \langle f \circ \eta \rangle \), so \( \overline{L}f \) is a digraph map. 

For morphisms \( f \colon G \to H \) and \( g \colon H \to K \), it is easy to verify that
\[
\overline{L}(g \circ f) = \overline{L}(g) \circ \overline{L}(f), \quad \text{and} \quad \overline{L}(\mathrm{id}_G) = \mathrm{id}_{\overline{L}G}.
\]
Thus, \( \overline{L} \colon \mathcal{DG}^* \to \mathcal{DG}^* \) is a functor.
\end{proof}
\end{proposition}
Next, we explore the relation between the digraph homotopy groups  $\overline{\pi}_n(G)$ of a based digraph $G$ and its loop-digraph $\overline{L}G$. At first, we construct a duality map in  Proposition \ref{duli}, which will induce an isomorphism $\Phi^{n+1}\colon \overline{\pi}_{n}(\overline{L}G)\rightarrow \overline{\pi}_{n+1}(G).$ 

\begin{definition}
The \emph{mapping digraph} \((H,B)^{(G,A)}\) is defined as the digraph whose vertices are the elements of \(\operatorname{Hom}((G,A),(H,B))\), and there is an arrow \(f \to g\) if, for all  \(v \in V(G)\), either \(f(v) \to g(v)\) or \(f(v) = g(v)\). In particular, we consider the case \((H, \emptyset)^{(G, \emptyset)}\), which we denote simply by \(H^G\).
 \end{definition}
\begin{proposition}\label{duli}
 There is an isomorphism $$\phi\colon  \operatorname{Hom}((G,A)\Box (G^{'},A^{'}); (H,B))\rightarrow \operatorname{Hom}((G,A);((H,B)^{(G^{'},A^{'})}, B^{G^{'}}))$$ defined by $f\mapsto \phi(f)$, where $\phi(f)(v)(v') = f(v,v')$.
 \begin{proof}
 First, for any $f\in Hom((G,A)\Box (G^{'},A^{'});(H,B))$, $\phi(f)\in Hom((G,A);((H,B)^{(G^{'},A^{'})}, B^{G^{'}}))$. Now let us check $\phi(f):G\rightarrow (H,B)^{(G^{'},A^{'})}$ is a digraph map. Let  $f_{g}=\phi(f)(g)$. We need to check $f_{g}:(G^{'},A^{'})\rightarrow (H,B)$ is a relative digraph map. Fixing $g$, observe that $f_{g}(g^{'}_{1}\rightarrow g^{'}_{2})$ is an arrow $f(g,g^{'}_{1})\rightarrow f(g,g^{'}_{2}) $ or a vertex $f(g,g^{'}_{1})$ by the definition of box product. Thus $f_{g}$ is a digraph map. It is easily to check that $f_{g}:(G^{'},A^{'})\rightarrow (H,B)$ is a relative digraph map.

 Next we check that if $g\rightarrow \widetilde{g}$, then $f_g\rightarrow f_{\widetilde{g}}$ or $f_g=f_{\widetilde{g}}$. Since for any $g^{'}\in G^{'}$, $f_g(g^{'})\rightarrow f_{\widetilde{g}}(g^{'})$ or $f_g(g^{'})= f_{\widetilde{g}}(g^{'})$, we obtain $f_g\rightarrow f_{\widetilde{g}}$. Thus $\phi(f)$ is a digraph map from $G$ to $(H,B)^{(G^{'},A^{'})}$. Clearly $\phi(f)(A)\subseteq B^{G^{'}}$. Hence $\phi(f)\in Hom((G,A);((H,B)^{(G^{'},A^{'})}, B^{G^{'}}))$.

Then we prove $\phi $ is surjective. For any $f\in Hom((G,A);((H,B)^{(G^{'},A^{'})}, B^{G^{'}}))$, $f_{g}:(G^{'},A^{'})\rightarrow (H,B)$ is a digraph map for any $g\in G$, and if $g\in A$, then $f_{g}:G^{'}\rightarrow B$. Then we define a map $F:(G,A)\Box (G^{'},A^{'})\rightarrow (H,B)$ by $(g,g^{'})\mapsto f_{g}(g^{'})$. Observe that $F\in Hom((G,A)\Box (G^{'},A^{'});(H,B))$ and $\phi(F)=f$.  Obviously, $\phi$ is injective. Hence $\phi$ is an isomorphism.
 \end{proof}
 \end{proposition}
By this proposition, there is an isomorphism $$\delta_{n+1}\colon \operatorname{Hom}((J_{m_i},\partial J_{m_i})^{\Box (n+1)};(G,\ast))\approx \operatorname{Hom}((J_{m_i},\partial J_{m_i})^{\Box n};((G,\ast)^{(J_{m_{n+1}},\partial J_{m_{n+1}})},\ast^{J_{m_{n+1}}}).$$
As $(G,\ast)^{(J_{m_{n+1}},\partial J_{m_{n+1}})}$ is a sub-digraph of $LG$, it induces a relative digraph map $$i_{m_{n+1}}\colon ((G,\ast)^{(J_{m_{n+1}},\partial J_{m_{n+1}})}, \ast^{J_{m_{n+1}}})\rightarrow (LG,\ast^{J_{m_{n+1}}})$$ defined by $f\mapsto f.$ Notice that there is a surjective digraph map
$$p_{m_{n+1}}\colon  (LG,\ast^{J_{m_{n+1}}})\rightarrow (\overline{L}G,\langle l_{\ast} \rangle)$$ such that $f\mapsto \langle f\rangle.$
Then there is a digraph map $$p_{m_{n+1}}\circ i_{m_{n+1}}\colon ((G,\ast)^{(J_{m_{n+1}},\partial J_{m_{n+1}})}, \ast^{J_{m_{n+1}}})\rightarrow (\overline{L}G,\langle l_{\ast} \rangle),$$
which induces a map $$\operatorname{Hom}((J_{m_i},\partial J_{m_i})^{\Box n};((G,\ast)^{(J_{m_{n+1}},\partial J_{m_{n+1}})}, \ast^{J_{m_{n+1}}}))\stackrel{P_{m_{n+1}}}{\rightarrow} \operatorname{Hom}((J_{m_i},\partial J_{m_i})^{\Box n};(\overline{L}G,\langle l_{\ast} \rangle)).$$
Combining this map with Proposition \ref{duli}, we have  $$ \operatorname{Hom}((J_{m_{i}},\partial J_{m_{i}})^{\Box (n+1)};(G,\ast)) \stackrel{\phi_{m_{n+1}}:= P_{m_{n+1}}\circ \phi}{\longrightarrow}
\operatorname{Hom}((J_{m_{i}},\partial J_{m_{i}})^{\Box n};(\overline{L}G,\langle l_\ast\rangle)).$$ By taking direct limit, $\phi_{m_{n+1}}$ induces a map $$\Phi^{n+1}\colon  \operatorname{Hom}((J,\partial J)^{\Box (n+1)};(G,\ast)) \rightarrow
\operatorname{Hom}((J,\partial J)^{\Box n};(\overline{L}G,\langle l_\ast\rangle)).$$

As a special case, if $n=0$, $\Phi^1 \colon  \operatorname{Hom}((J,\partial J); (G,\ast))\rightarrow \operatorname{Hom}(0;\overline{L}G)$ should be defined by $\gamma\mapsto e_{\gamma}$, where $e_{\gamma}\colon 0\rightarrow \overline{L}G$ such that $e_\gamma(0) =\langle\gamma\rangle$ since $(J_m,\partial J_m) = (0,\emptyset)$ $\Box(J_m,\partial J_m)$.

Next, along the idea of adjusting the length of representative in subdivision class in \cite[Proposition 7.4]{A-1}, we obtain the following lemma by using our common subdivision in Lemma \ref{subre}.
\begin{lemma}\label{lattice}
Let $G$ be a based digraph with base-point $\ast$. For any digraph map $\widetilde{f}\colon (J_{m_i},\partial J_{m_i})^{\Box n}\rightarrow (\overline{L}G,\langle l_\ast\rangle)$, there is an integer $m_{n+1}$ and a digraph map 
$$f\colon (J_{m_i},\partial J_{m_i})^{\Box (n+1)}\rightarrow (G,\ast)$$ such that $\phi_{m_{n+1}}(f)=\widetilde{f}$ for $n\geq1$.
\begin{proof}
Denote $\widetilde{f}(a)$ by $\langle\widetilde{f}_a\rangle$ for each $a\in (J_{m_{i}},\partial J_{m_{i}})^{\Box n}$, where we refer to $$\widetilde{f}_a\colon (J_{m_{n+1}},\partial J_{m_{n+1}})\rightarrow (G,\ast)$$ a standard loop by Lemma~\ref{standard}.
The lemma holds by setting 
$$f:(J_{m_1},\partial J_{m_1})^{\Box n}\Box (J_{m_{n+1}},\partial J_{m_{n+1}})\rightarrow (G,\ast)$$
with $f(a,b)=\widetilde{f}_a(b)$ for $$(a,b)\in V((J_{m_i},\partial J_{m_i})^{\Box n}\Box (J_{m_{n+1}},\partial J_{m_{n+1}}))= V((J_{m_i},\partial J_{m_i})^{\Box n})\times V (J_{m_{n+1}},\partial J_{m_{n+1}})).$$
\end{proof}
\end{lemma}
\begin{remark}\label{plattice}
For any digraph map
$\widetilde{f}\colon (J_{m_i},\partial J_{m_i})^{\Box n}\rightarrow (\overline{P}G,\langle l_\ast\rangle),$ there is still a digraph map
$$f\colon (J_{m_i},\partial J_{m_i})^{\Box n}\Box (J_{m_{n+1}},0)\rightarrow (G,\ast)$$
such that $p'_{m_{n+1}}\circ j_{m_{n+1}}\circ \phi(f) =\widetilde{f},$ where
$$j_{m_{n+1}}\colon (G^{(J_{m_{n+1}},0)},\ast^{J_{m_{n+1}}})\rightarrow (PG,\ast^{J_{m_{n+1}}})$$ is an
embedding digraph map and $$p'_{m_{n+1}}\colon (PG,\ast^{J_{m_{n+1}}})\rightarrow(\overline{P}G,\langle l_{\ast}\rangle)$$ maps $f$ to $\langle f\rangle.$
\end{remark}
Inspired by Lemma \ref{lattice}, we explore the relationship between $\overline{\pi}_{n+1}(G)$ and $\overline{\pi}_n(\overline{L}G)$ further.
\begin{proposition}\label{main}
The map $$\Phi^{n+1}\colon  \operatorname{Hom}((J,\partial J)^{\Box (n+1)};(G,\ast)) \rightarrow
\operatorname{Hom}((J,\partial J)^{\Box n};(\overline{L}G,\langle l_\ast\rangle)) $$ induces an isomorphism $$\Phi^{n+1}_{\ast}\colon \overline{\pi}_{n+1}(G) \rightarrow \overline{\pi}_n(\overline{L}G),\text{ }[f]\mapsto [\Phi^{n+1}(f)]$$ for $n\geq 1.$
\begin{proof}
First, we check that $\Phi^{n+1}_{\ast}$ is well-defined. Let $$f\colon (J_{m_i},\partial J_{m_i})^{\Box (n+1)}\rightarrow (G,\ast)$$ and $$g\colon (J_{l_i},\partial J_{l_i})^{\Box (n+1)}\rightarrow (G,\ast)$$ be representatives in $\{f\}$, $\{g\}\in \operatorname{Hom}((J,\partial J)^{\Box (n+1)};(G,\ast)).$ Recall $\{f\}\s_{1} \{g\}$ or $\{f\}\s_{-1} \{g\}$ if and only if $f\s_1 g$ or $f\s_{-1}g$ respectively for elements $\{f\},\text{ }\{g\}$ $\in \operatorname{Hom}((J,\partial J)^{\Box (n+1)};(G,\ast))$ regardless of the representatives.
%So we do not distinguish between $f$ and $\{f\}$ in what follows.
Suppose  $f\simeq_{F}g$, then there is a sequence of maps $\{f_{l}\}_{l=0}^{s}$ such that $f_{l}\s_1 f_{l+1}$ or $f_{l}\s_{-1} f_{l+1}$. We need to verify $ \phi_{m_{n+1}}(f)\simeq_{F}\phi_{l_{n+1}}(g)$. Then $\Phi_{{n+1}}(f)\simeq_{F}\Phi_{{n+1}}(g).$ It is sufficient to show that if $f\simeq_{1}g$ or $f\simeq_{-1}g$, then $ \phi_{m_{n+1}}(f)\simeq_{F}\phi_{l_{n+1}}(g)$. Denote   $\widetilde{f}=\phi_{m_{n+1}}(f)$ and $\widetilde{g}=\phi_{l_{n+1}}(g).$

If $f\simeq_{1}g$, then there exist two subdivisions $$\overline{f}\colon (J_{M_i},\partial J_{M_i})^{\Box (n+1)}\rightarrow (G,\ast)$$ and $$\overline{g}\colon (J_{M_i},\partial J_{M_i})^{\Box (n+1)}\rightarrow (G,\ast)$$ of $f$ and $g$ respectively such that $\overline{f}\dr \overline{g} $, thus $\phi_{M_{n+1}}(\overline{f})=\widetilde{\overline{f}}\dr \widetilde{\overline{g}}=\phi_{M_{n+1}}(\overline{g})$. Now let us consider the relationship between $\widetilde{\overline{g}}$ and $\widetilde{g}$ in the following three cases.

\begin{enumerate}
  \item If $\overline{g}$ is a subdivision of $g$ in the first $n$ coordinates, then $\widetilde{\overline{g}}$ is a subdivision of $\widetilde{g}$. Therefore $\widetilde{\overline{g}}\s_1 \widetilde{g}.$
  \item If $\overline{g}$ is a subdivision of $g$ in the last coordinate, then $\widetilde{\overline{g}}_{a}$ is a subdivision of $\widetilde{g}_a$ for all $a\in (J_{M_{i}},\partial J_{M_{i}})^{\Box n}$, so $\langle\widetilde{\overline{g}}_{a}\rangle = \langle\widetilde{g}_a\rangle$ in $\overline{L}G$. Hence $\widetilde{\overline{g}} = \widetilde{g}$.
  \item By the subdivision decomposition, if $\overline{g}$ is a subdivision of $g$ not only in the first $n$ coordinates but also in the last, then $\widetilde{\overline{g}}\s_1 \widetilde{g}$ by considering both cases above.
\end{enumerate}

In conclusion, $\widetilde{\overline{g}}\s_1 \widetilde{g}$ or $\widetilde{\overline{g}}= \widetilde{g}$.  Similarly,  $\widetilde{\overline{f}}\s_1\widetilde{f}$ or $\widetilde{\overline{f}}=\widetilde{f}$. Therefore $\widetilde{f}\s_F \widetilde{g}.$ Hence if $f\simeq_{F}g$ then $\widetilde{f}\simeq_{F}\widetilde{g}$. That is to say, $\Phi^{n+1}_{\ast}$ is well-defined.

Since for any $\widetilde{f}\colon  (J_{m_{i}},\partial J_{m_{i}})^{\Box n}\rightarrow (\overline{L}G,\langle l_\ast\rangle)$,
there is a digraph map $f\colon (J_{m_{i}},$ $\partial J_{m_{i}})^{\Box (n+1)}\rightarrow (G,\ast)$ such that $\phi_{m_{n+1}}(f)=\widetilde{f}$ by Lemma \ref{lattice}. Hence $\Phi^{n+1}$ is surjective, so does $\Phi^{n+1}_{\ast }$.

Next let us show $\Phi^{n+1}_{\ast }$ is injective. We only need to show that $\{f\}\simeq_{F}\{g\}$ if $[\Phi^{n+1}(\{f\})]=[\Phi^{n+1}(\{g\})]$. Since $\Phi^{n+1}(\{f\}) = \{\phi_{m_{n+1}}(f)\}$ and $\Phi^{n+1}(\{g\}) = \{\phi_{l_{n+1}}(g)\}$,
if $\Phi^{n+1}(\{f\})\simeq_{F}\Phi^{n+1}(\{g\})$, then there is a sequence of digraph maps $\{\widetilde{f}_i\}_{i=0}^{l}$ from $\phi_{m_{n+1}}(f)$ to $\phi_{l_{n+1}}(g)$ such that $\widetilde{f}_{i}\s_1\widetilde{f}_{i+1}$ or $\widetilde{f}_{i}\s_{-1}\widetilde{f}_{i+1}$ for $0\leq i\leq l-1$.

Suppose $$\widetilde{f}_0\colon (J_{m_{i}},\partial J_{m_i})^{\Box n}\rightarrow (\overline{L}G,\langle l_\ast\rangle)$$
and $$\widetilde{f}_1\colon (J_{k_{i}},\partial J_{k_i})^{\Box n}\rightarrow (\overline{L}G,\langle l_\ast\rangle),$$
then there exist $$f_0\colon (J_{m_{i}},\partial J_{m_i})^{\Box (n+1)}\rightarrow (G,\ast)$$ and
$$f_1\colon (J_{k_{i}},\partial J_{k_i})^{\Box (n+1)}\rightarrow (G,\ast)$$ such that
$\phi_{m_{n+1}}(f_0) = \widetilde{f}_0$ and $\phi_{k_{n+1}}(f_1) = \widetilde{f}_1$ by Lemma
\ref{lattice}. To prove $\Phi^{n+1}_{\ast}$ is injective, we only need to prove $f_0\s_1 f_1$ if
$\widetilde{f}_{0}\s_1 \widetilde{f}_1$ and $f_0\s_{-1} f_1$ if $\widetilde{f}_{0}\s_{-1} \widetilde{f}_1$. By definition of one-step $F$-homotopy, there are subdivisions
$\overline{\widetilde{f}}_0$ and $\overline{\widetilde{f}}_1$ of $\widetilde{f}_{0}$ and
$ \widetilde{f}_1$ by $q_0$ and $q_1$ respectively such that
$\overline{\widetilde{f}}_0 \dr \overline{\widetilde{f}}_1$,
where $\overline{\widetilde{f}}_0, ~\overline{\widetilde{f}}_1\colon (J_{M_{i}},\partial J_{M_i})^{\Box n}\rightarrow (\overline{L}G, \langle l_\ast\rangle )$.
Thus there exist $$\overline{f}_0:= f_0\circ(q_0\Box id)\colon (J_{M_{i}},\partial J_{M_i})^{\Box n}\Box (J_{m_{n+1}},\partial J_{m_{n+1}})\rightarrow (G,\ast)$$
 and $$\overline{f}_1:=f_1\circ(q_1\Box id)\colon (J_{M_{i}},\partial J_{M_i})^{\Box n}\Box (J_{k_{n+1}},\partial J_{k_{n+1}})\rightarrow (G,\ast)$$
such that $\phi_{m_{n+1}}(\overline{f}_0) = \overline{\widetilde{f}}_0$ and $\phi_{k_{n+1}}(\overline{f}_1) = \overline{\widetilde{f}}_1$. It is clear that $\overline{f}_0$ and $\overline{f}_1$ are subdivisions of $f_0$ and $f_1$ respectively.

Since $\overline{\widetilde{f}}_0 \dr \overline{\widetilde{f}}_1$, there is a digraph map $$F\colon  (J_{M_{i}},\partial J_{M_i})^{\Box n}\Box J_{1}\Box (J_{M_{n+2}},\partial J_{M_{n+2}})\rightarrow (G,\ast) $$ such that $$\phi_{M_{n+2}}(F)|_{(J_{M_{i}},\partial J_{M_i})^{\Box n}\Box 0} = \overline{\widetilde{f}}_0$$ and $$\phi_{M_{n+2}}(F)|_{(J_{M_{i}},\partial J_{M_i})^{\Box n}\Box 1} = \overline{\widetilde{f}}_1$$ by Lemma \ref{lattice}.
Denoting $$\overline{f}_{0}':=F|_{(J_{M_{i}},\partial J_{M_i})^{\Box n}\Box 0\Box (J_{M_{n+2}},\partial J_{M_{n+2}})}$$ and $$\overline{f}_1':=F|_{(J_{M_{i}},\partial J_{M_i})^{\Box n}\Box 1\Box (J_{M_{n+2}},\partial J_{M_{n+2}})},$$ we obtain $\overline{f}_0'\dr \overline{f}_1'$.

Next we claim $\overline{f}_0'$ is a subdivision of $\overline{f}_0$. For any $a\in (J_{M_{i}},\partial J_{M_i})^{\Box n}$,
$(\overline{f}_{0})_a\colon  (J_{m_{n+1}},$ $\partial J_{m_{n+1}})\rightarrow(G,\ast) $ and $(\overline{f}'_{0})_a\colon  (J_{M_{n+2}},\partial J_{M_{n+2}})\rightarrow(G,\ast) $ such that $\langle(\overline{f}'_{0})_a\rangle= \langle(\overline{f}_{0})_a\rangle$, there is a common subdivision by some shrinking map $q\colon (J_{M_{n+2}},\partial J_{M_{n+2}})\rightarrow (J_{m_{n+1}},$ $\partial J_{m_{n+1}})$ for all $a$. Hence $\overline{f}'_0$ is the subdivision of $\overline{f}_0$ by $id^{\Box n}\Box q$. As  $\overline{f}_0$ is a subdivision of $f_0$, it follows that $\overline{f}'_0$ is a subdivision of $f_0$. Similarly, $\overline{f}'_1$ is a subdivision of $f_1$. Thus $f_0\s_1 f_1.$ In a similar way, $f_0\s_{-1} f_1$ if $\widetilde{f}_{0}\s_{-1} \widetilde{f}_1$. It follows that if $\widetilde{f}_{0} = \phi_{m_{n+1}}(f)\simeq_{F}\phi_{l_{n+1}}(g) = \widetilde{f}_{l}$, then $f\simeq_{F}g$. Hence $\Phi^{n+1}_{\ast }$ is injective.

Finally, we show that $\Phi^{n+1}_{\ast}$ is a homomorphism. It is sufficient to show $[\{\phi_{M_{n+1}}(f\cdot g)\}]= [\{\phi_{m_{n+1}}(f)\}\cdot\{\phi_{l_{n+1}}(g)\}],$ where $M_{n+1} = max\{m_{n+1},~l_{n+1}\}$.
By the definition of the multiplication, $ f\cdot g =f'\vee g'$, that is, 
$$ f\cdot g \colon  (J_{(m_{1}+l_1)},\partial J_{(m_{1}+l_1)})\Box\cdots \Box (J_{M_j},\partial J_{M_j})\Box\cdots \Box(J_{M_{n+1}},\partial J_{M_{n+1}}) \rightarrow (G,\ast),$$
where $f^{'}$ and $g^{'}$ are the subdivisions of $f$ and $g$ respectively such that $f'\vee g'$ is still a grid map and $M_i=max\{m_i,~l_i\}$ for $2\leq i\leq n+1$.
Then
$$  \phi_{M_{n+1}}(f\cdot g)(i_1,i_2,...,i_{n})= \left\{
                     \begin{array}{ll}
                      \phi_{M_{n+1}}(f^{'})(i_1,i_2,...,i_n), & \hbox{$i_1\leq m_1$;} \\
                        \phi_{M_{n+1}}(g^{'})(i_1-m_1,i_2,...,i_n), & \hbox{$i_1> m_1$;}
                     \end{array}
                   \right
.$$
and
$$  (\phi_{m_{n+1}}(f)\cdot \phi_{l_{n+1}}(g))(i_1,i_2,...,i_{n})= \left\{
                     \begin{array}{ll}
                       (\phi_{M_{n+1}}(f))^{'}(i_1,i_2,...,i_n), & \hbox{$i_1\leq m_1$;} \\
                       (\phi_{M_{n+1}}(g))^{'}(i_1-m_1,i_2,...,i_n), & \hbox{$i_1> m_1$.}
                     \end{array}
                   \right
.$$
Since for any $(i_1,i_2,...,i_n)\in (J_{(m_{1}+l_1)},\partial J_{(m_{1}+l_1)})\Box (J_{M_2},\partial J_{M_2})\Box\cdots \Box(J_{M_{n}},\partial J_{M_{n}})$, $$\phi_{M_{n+1}}(f\cdot g)(i_1,i_2,...,i_{n})=(\phi_{m_{n+1}}(f)\cdot \phi_{l_{n+1}}(g))(i_1,i_2,...,i_{n}).$$ Then $\{\phi_{M_{n+1}}(f\cdot g)\}= \{\phi_{m_{n+1}}(f)\}\cdot\{\phi_{l_{n+1}}(g)\}.$ Thus $\Phi^{n+1}(f\cdot g) = \Phi^{n+1}(f)\cdot \Phi^{n+1}(g),$ so $\Phi^{n+1}_{\ast}$ is a homomorphism. Putting it all together, $\Phi^{n+1}_{\ast}$ is an isomorphism.
\end{proof}
\end{proposition}
Due to the fact that $\Phi^{n+1}_{\ast}$ is an isomorphism and the argument in the previous proof, we can precisely describe the inverse homomorphism $\Psi^{n+1}_{\ast}$ of $\Phi^{n+1}_{\ast}$ for $n\geq 1.$
\begin{corollary}The inverse $\Psi^{n+1}_{\ast}$ of $\Phi^{n+1}_{\ast}$ is the map
$$\Psi^{n+1}_{\ast}\colon  \overline{\pi}_n(\overline{L}G) \rightarrow \overline{\pi}_{n+1}(G), \text{ }[\widetilde{f}] \mapsto [f],$$ where $f$ is the digraph map constructed by Lemma \ref{lattice} for $n\geq1$. \qed
\end{corollary}
Denote $\overline{L}^{n}G= \overline{L}(\overline{L}^{n-1}G)$. By Proposition \ref{main}, we have the following theorem.
\begin{theorem}\label{dulo}Let $G$ be a based digraph. Then
$\overline{\pi}_{1}(\overline{L}^{n}G) \approx \overline{\pi}_{n+1}(G)$ for all $n\geq 1$.
\begin{proof}
Clearly, $\overline{\pi}_2(G)\approx \pi_1(\overline{L}G)$. Assume $\overline{\pi}_{n}(G)\approx \overline{\pi}_1(\overline{L}^{n-1}G)$, then due to the fact $\overline{\pi}_{n+1}(G)\approx \overline{\pi}_{n}(\overline{L}G),$ we obtain the result.
\end{proof}
\end{theorem}
In summary, we have constructed a reduced loop-digraph $\overline{L}G$ for a based digraph $G$ such that
$\overline{\pi}_{1}(\overline{L}^{n}G) \approx \overline{\pi}_{n+1}(G)$.

To this point, we have obtained
some similar properties to topological spaces. But another natural questions appears, ``Is our definition of  homotopy groups useful? Does it indeed distinguish different digraphs?". Let us see the example below to answer the above questions.
\begin{example}\label{sphere}
 Let $G$ be the based digraph shown as follows. We want to prove $\overline{\pi}_2(G)\neq 0$. To do so, we only need to find a non-trivial element $[f]$ in $\overline{\pi}_2(G)$.
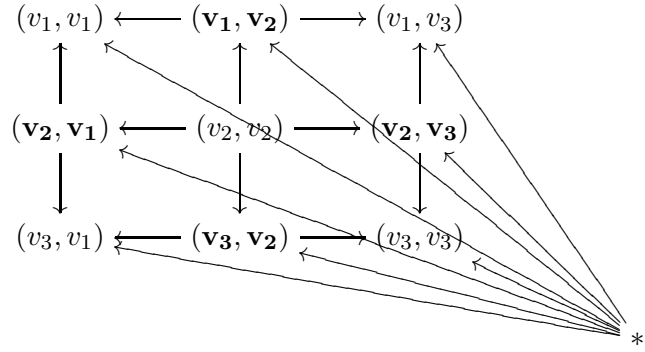
$$
\xymatrix@=2em{
  (v_1,v_1)  & \mathbf{(v_1,v_2)} \ar[l] \ar[r] & (v_1,v_3) & & \\
   \mathbf{(v_2,v_1)}\ar[u]\ar[d]  & (v_2,v_2) \ar[d]\ar[l]\ar[u] \ar[r] & \mathbf{(v_2,v_3)} \ar[u]\ar[d] & &\\
  (v_3,v_1)  & \mathbf{(v_3,v_2)} \ar[r]\ar[l] & (v_3,v_3)  &  & \\
  &                          &   &  & \ast\ar[ull]\ar[ulll]\ar[ullll]\ar[uull]\ar[uuulll]\ar[uuull]\ar[uuullll]\ar[uullll] }
  $$
Define $$f\colon (J_{4},\partial J_4)\Box (J_4,\partial J_4)\rightarrow (G,\ast)$$ by $f(i,j)=(v_i,v_j)$ when $1\leq i,\text{ }j\leq 3$. We claim $f$ is not null-homotopic.

Suppose $f\s_F e_{\ast}$, then there is a family of digraph maps $\{f_i\}_{i=0}^{l}$ such that $f_0 = f$, $f_l=e_{\ast}$ and
$f_i\s_1 f_{i+1}$ or $f_i\s_{-1}f_{i+1}$ for $1\leq i\leq l-1,$
where $e_{\ast}\colon J_4^{\Box 2}\rightarrow \ast.$ Then $f_{l-1}$ must be a digraph map such that
 $f_{l-1}\s_{-1} e_{\ast}.$ Each digraph map $f_i$ can be seen as a sub-digraph $f_i(G)$
of $G$, therefore $\{f_i\}_{i=0}^{l}$ is a sequence of sub-digraphs of $G$.  Since the   vertex
$(v_2,v_2)\in G$ does not connect with the base-point $\ast,$ the center vertex in the sub-digraph $f_{l-1}(G)$ should
be bold or $\ast$. That is to say, there exist two adjacent sub-digraphs $f_s(G)$ and $f_{s+1}(G)$ in $\{f_i(G)\}_{i=0}^l$ turning the   vertex $(v_2,v_2)$ into some bold vertex.

Since the  bold vertices lie in symmetric positions, we assume that the   vertex $(v_2,v_2)$ is turned into the  bold vertex $(v_1,v_2)$ in passing from $f_s$ to $f_{s+1}$.
Here we claim that the center vertex $(v_2,v_2)$ and the vertices connecting the vertex $(v_2,v_2)$ in $f_s(G)$ are the same as $f_0(G)$. Then to guarantee that $f_{s+1}$ is a digraph map, the vertices connecting with the vertex $(v_1,v_2)$ in $f_{s+1}(G)$ should be different from $f_s(G)$ since $(v_1,v_2)$ is not connected with the other three  bold vertices. That is to say, some  bold vertex is changed in passing from $f_{s}(G)$ to $f_{s+1}(G)$. Suppose the  bold vertex $(v_3,v_2)$ is changed in passing from $f_s(G)$ to $f_{s+1}(G)$, then the vertex $(v_3,v_2)$ can change to $(v_3,v_1)$ or $(v_3,v_3)$, which are not connect with the center vertex $(v_1,v_2)$ in $f_{s+1}(G)$. So if there exists a one-step $F$-homotopy from $f_s$ to $f_{s+1}$, then the vertex $(v_1,v_2)$ connects with the vertex $(v_3,v_1)$ or $(v_3,v_3)$. This is a contradiction. Hence the   vertex $(v_2,v_2)$ can not change in passing from $f_s(G)$ to $f_{s+1}(G)$, that is to say, there is no integer $s$ such that $f_s(G)$ sends the   vertex $(v_2,v_2)$ to some  bold vertex. Thus there does not exist a digraph map $f_{l-1}$ such that $f_{l-1}\s_{-1}e_{\ast}$. Hence $f$ is not null-homotopic.

Now let us prove the claim that the center vertex $(v_2,v_2)$ and the vertices connecting the vertex $(v_2,v_2)$ in $f_s(G)$ are the same as $f_0(G)$. It is sufficient to prove if the   vertex $(v_2,v_2)$ does not change, then the  bold vertices don't change.
If we prove that between the adjacent digraph maps $f_i$ and $f_{i+1}$ the  bold vertices don't change when the   vertex $(v_2,v_2)$ does not change, then $f_s(G)$ is a sub-digraph preserving center and  bold vertices. Suppose the   vertex $(v_2,v_2)$ is preserved between the adjacent sub-digraphs, then each  bold vertex can only change into a  bold vertex since $(v_2,v_2)$ only connects with the  bold vertices. But any two  bold vertices do not connect with each other. This is a contradiction. So the center vertex $(v_2,v_2)$ and the vertices connecting the vertex $(v_2,v_2)$ in $f_s(G)$ are same the as $f_0(G)$.
\end{example}
Observe the digraph $G$ in Example \ref{sphere} is similiar to the topological space $S^2$. However, this may be a misleading comparison. The homotopy groups of $G$ are mysterious and need to be investigated further.

 Now let us see the relation between the loop-digraph $LG$ and the reduced loop-digraph $\overline{L}G$ of $G$.
\begin{proposition}\label{0i}
There is a based digraph map $$p\colon  (LG,l_{\ast})\rightarrow (\overline{L}G,\langle l_{\ast}\rangle),\text{ } f\mapsto \langle f\rangle,$$
which induces a surjective morphism $$p_{n}\colon  \overline{\pi}_{n}(LG)\rightarrow \overline{\pi}_{n}(\overline{L}G)$$ for $n\geq 0$. In the special case when $n=0,$ $p_0$ is an isomorphism.
\begin{proof}
Obviously, if $f_1\rightarrow f_2$ in $LG$, then $\langle f_1\rangle\rightarrow \langle f_2\rangle$ or $\langle f_1\rangle =  \langle f_2\rangle$ by definition of $\overline{L}G$. Hence $p$ is a digraph map preserving the base-point. Therefore $p$ induces a morphism $p_n\colon \overline{\pi}_{n}(LG)\rightarrow \overline{\pi}_{n}(\overline{L}G)$ by Proposition \ref{map}. Clearly $p_n$ is surjective.

Next we check $p_0$ is injective.
For any $[f],~[g]\in \overline{\pi}_{0}(LG)$, if $[p\circ f]= p_{0}([f])= p_{0}([g])$ $=[p\circ g]$, then there exists a line digraph $I_l$ and a digraph map $F\colon I_l\rightarrow \overline{L}G$ such that $F(0) = (p\circ f)(0)$ and $F(l) = (p\circ g)(0)$ in $\overline{L}G$. Let $F(i)= \langle\widetilde{f}_i\rangle$ for $0\leq i\leq l.$ If $\langle\widetilde{f}_0\rangle\rightarrow \langle\widetilde{f}_1\rangle$, then there exist two subdivisions $\overline{\widetilde{f}}_0$ and $\overline{\widetilde{f}}_1$ of $\widetilde{f}_0$ and $\widetilde{f}_1$ respectively such that $\overline{\widetilde{f}}_0 \dr \overline{\widetilde{f}}_1$. Similarly for $\langle\widetilde{f}_1\rangle\leftarrow \langle\widetilde{f}_2\rangle$, we also have subdivisions $\overline{\widetilde{f}}^{'}_1$ and $\overline{\widetilde{f}}_2$ of $\widetilde{f}_1$ and $\widetilde{f}_2$ respectively such that $\overline{\widetilde{f}}^{'}_1 \dl \overline{\widetilde{f}}_2$. Repeating this procedure, we obtain the diagram shown as follows.
\begin{small}$$\xymatrix@C=1em@R=0.7em{
                & \overline{\widetilde{f}}_0  \ar[dl] \ar@<.5ex>[r] \ar@<-0.5ex>[r]&    \overline{\widetilde{f}}_1 \ar[dr] & & \overline{\widetilde{f}}^{'}_1\ar[dl]& \overline{\widetilde{f}}_2\ar@<.5ex>[l] \ar@<-0.5ex>[l]\ar[dr] & & ...& & \overline{\widetilde{f}}^{'}_{n-1}\ar@<.5ex>[r] \ar@<-0.5ex>[r]\ar[dl]& \overline{\widetilde{f}}_n\ar[dr] & \\
    f(0) = \widetilde{f}_0            &               &      &  \widetilde{f}_1 & & & \widetilde{f}_2 & ...& \widetilde{f}_{n-1}& &  &\widetilde{f}_n =g(0) }
$$
\end{small}
Therefore there exists a line digraph $I_{N}$ and a digraph map $\widetilde{F}\colon I_{N}\rightarrow LG$ such that $\widetilde{F}(0)=f(0)$ and $\widetilde{F}(N)=g(0)$ in $LG$. Thus $[f]=[g]$ in $\overline{\pi}_{0}(LG).$

\end{proof}
\end{proposition}
Proposition \ref{0i} immediately implies the following corollary.
\begin{corollary}
Let $G$ be a based digraph. Then the following hold: 

(1) $\pi_1(G)=\overline{\pi}_1(G)$;

(2) $p_1\colon  \pi_2(G)\rightarrow \overline{\pi}_2(G)$ is surjective.
\begin{proof}
(1) If $n=0$, then $\pi_1(G)=\pi_0(LG) = \overline{\pi}_0(LG) \stackrel{p_0}{\approx} \overline{\pi}_0(\overline{L}G) =\overline{\pi}_1(G),$ implying that $\pi_1(G)\approx\overline{\pi}_1(G)$.
Moreover as Lemma \ref{equ} says that $C$-homotopy is same as $F$-homotopy between paths,  we directly obtain the $\pi_1(G) =\bar{\pi}_1(G).$

 (2) If $n=1$,  then $\pi_2(G)=\pi_1(LG) \approx \overline{\pi}_1(LG)\stackrel{p_1}{\twoheadrightarrow} \overline{\pi}_1(\overline{L}G) =\overline{\pi}_2(G).$
\end{proof}
\end{corollary}
For more precise relation between $\pi_{n}(G)$ and $\overline{\pi}_{n}(G)$ for $n\geq 1$ is not clear. In the future, we will continue to explore this.
\section{Puppe Sequence of Digraphs}
For any continuous map $f$, one can construct a fibration that induces a long exact
sequence of homotopy groups, which is called the Puppe sequence \cite{Puppe} or fiber sequence \cite{Swiczer}.
This Puppe sequence is highly connected to fibration theory and fibre bundle theory in topological spaces.
 Inspired by this fact,
this section is concerned with establishing a Puppe sequence for based digraphs in a constructing method on underlying digraphs.
Following the idea for a Puppe sequence of spaces~\cite{Maunder}, we start with the contractibility of $\overline{P}G$.
To this point we can only prove $\overline{P}G$ is weakly contractible,
 that is, $\overline{\pi}_n(\overline{P}G)=0$
for all $n\geq 0$, instead of contractible.
\begin{proposition}\label{weak}
For any based digraph $G$, $\overline{\pi}_n(\overline{P}G)=0$ for $n\geq0$.
\begin{proof}
In particular, $\overline{\pi}_0(\overline{P}G)=0.$ To prove this, we need to show that for any digraph map $l\colon 0\rightarrow \overline{P}G$, $[l]=[\ast]$, where $\ast\colon  0\rightarrow \langle l_\ast\rangle$ and $l_{\ast}\colon  (J_{2},\partial J_2)\rightarrow (\ast,\ast)$. That is to say, there is a digraph map $F\colon  I_n\rightarrow \overline{P}G$ such that $F(0)= l(0)$ and $F(n) = \ast(0)=\langle l_\ast\rangle$.

Assume $l(0)=\langle\gamma\rangle$ and $\gamma\colon  (J_m,0)\rightarrow (G,\ast)$. We can construct a sequence of paths $\{\gamma_i\}_{i=1}^m$, where $\gamma_i = \gamma|_{J_{i}}\colon (J_{i},0)\rightarrow (G,\ast)$ for $1\leq i\leq m$. Clearly $\langle\gamma_i\rangle = \langle\gamma_{i+1}\rangle$ or $\langle\gamma_i\rangle $ $ \rightarrow \langle\gamma_{i+1}\rangle$ or $\langle\gamma_i\rangle \leftarrow \langle\gamma_{i+1}\rangle$. By the definition of $\gamma_i$,  $\gamma_m = \gamma$ and $\gamma_1 \s_{-1} l_{\ast}$. Therefore we can define a digraph map $F\colon J_{m}\rightarrow \overline{P}G$ by $F(i)= \langle\gamma_i \rangle$ such that $F(m)=\langle\gamma\rangle$ and $F(0)= \langle l_\ast\rangle $. Hence $[l] = [\ast]$. Since $l$ is arbitrary, $\overline{\pi}_0(\overline{P}G)=0$.

Now we turn to check $\overline{\pi}_{n}(G) = 0$ for $n\geq1$. It is sufficient to show that for any digraph map $l\colon (J_{m_i},\partial J_{m_i})^{\Box n}\rightarrow (\overline{P}G,\langle l_\ast\rangle)$, $l\s_F \widetilde{l}_{\ast}$, where $\widetilde{l}_{\ast}\colon  (J_{m_i},\partial J_{m_i})^{\Box n}\rightarrow (\langle l_\ast\rangle,$ $\langle l_\ast\rangle)$.

As in the proof of Lemma \ref{lattice}, one can easily check that for any $l\colon (J_{m_i},\partial J_{m_i})^{\Box n}\rightarrow (\overline{P}G,\langle l_\ast\rangle)$, there is a digraph map $$l'\colon (J_{m_i},\partial J_{m_i})^{\Box n}\Box (J_{m_{n+1}},0)\rightarrow (G,\ast)$$ such that $p\circ \phi(l') = l$, where $p\colon (PG, \ast^{J_{m_{n+1}}})\rightarrow (\overline{P}G, \langle l_{\ast}\rangle)$ sends $\gamma$ to $\langle \gamma\rangle.$ Then we can construct a map $$F\colon (J_{m_i},\partial J_{m_i})^{\Box n}\Box J_{m_{n+1}}\rightarrow (\overline{P}G,\langle l_{\ast}\rangle)$$ such that $F_{(J_{m_i},\partial J_{m_i})^{\Box n}\Box j} = f_j$, $1\leq j\leq m_{n+1}$, where $f_j :=p\circ \phi(l'|_{ (J_{m_i},\partial J_{m_i})^{\Box n} \Box (J_{j},0)})\colon $ $ (J_{m_i},\partial J_{m_i})^{\Box n} $ $\rightarrow (\overline{P}G,\langle l_\ast\rangle)$. Clearly $f_j\dr f_{j+1}$ or $f_j\dl f_{j+1}$, so $F$ is a digraph map such that \\
$F|_{(J_{m_i},\partial J_{m_i})^{\Box n}\Box m_{n+1}} = l$ and $F|_{(J_{m_i},\partial J_{m_i})^{\Box n}\Box 1} \dl \widetilde{l}_{\ast}$. Hence, $l\s_F \widetilde{l}_{\ast}$. Due to $l$ being arbitrary, $\overline{\pi}_n(\overline{P}G)=0$ for $n\geq1$. To sum up, $\overline{\pi}_n(\overline{P}G)$ is trivial for all $n\geq 0.$
\end{proof}
\end{proposition}

To build the Puppe sequence of digraphs, we should introduce the mapping path-digraph of a based digraph $G$. Some notations are required. For any based digraph $G$, there is a digraph map $e\colon \overline{P}G \rightarrow G$ defined by $\langle\lambda\rangle\mapsto \lambda(n)$, where $\lambda\colon  (J_n, 0)\rightarrow (G,\ast)$. The digraph map $e$ is called the \emph{evaluation map} of $G$. Since any two different subdivisions of the path $\lambda\colon (J_n, 0)\rightarrow (G,\ast)$ have the same end vertex, the map $e$ is well-defined. In what follows, we will write $\lambda(end)$, instead of $e(\langle \lambda\rangle)$. By the definition of the reduced path-digraph $\overline{P}G$, the map $e$ is a digraph map.

Let $G= (V_G,A_G)$ and $H=(V_H, A_H)$. The \emph{Cartesian product} $G\times H$ of $G$ and $H$ is the digraph whose vertex set
is $V_G\times V_H$ and whose arrow set contains the arrow $(v,w)\rightarrow (v',w')$ if and only if
\begin{itemize}
  \item $v=v'$ and $w\rightarrow w'$, or
  \item $v\rightarrow v'$ and $w=w'$, or
  \item $v\rightarrow v'$ and $w\rightarrow w'$.
\end{itemize}
Let us see the precise definition of pullback in the digraph category and check the universal property.
%Similar to topological spaces, we call a digraph $G$ \emph{path-connected} if for any two vertices $a$ and $b$ in $G$, there exists a path $\gamma: I_n\rightarrow G$ such that $\gamma(0) = a$ and $\gamma(n) =b$.
%\begin{prop}
%Given a digraph $G$, if we choose two different vertices $a$ and $b$ in one path component as base-point, then $\overline{P}G^{a}$ is homotopy equivalent to $\overline{P}G^{b}.$
%\end{prop}
\begin{definition}
Let $f\colon X\rightarrow Z$ and $g\colon  Y\rightarrow Z $ be two digraph maps. We call the triple $(X\times_{Z} Y, p_1, p_2 )$ in the following diagram the pullback of $f$ and $g$,
$$
\xymatrix{
  X\times_{Z} Y \ar[d]_{p_1} \ar[r]^{p_2}\pushoutcorner
                & Y \ar[d]^{g}  \\
  ~~~X~~~\ar[r]^{f}
               & Z             }
$$
 where $X\times_{Z} Y$ is the induced sub-digraph of the Cartesian product $X\times Y$ generated by the vertex set $V(X\times_{Z} Y)=\{(x,y)\in V(X)\times V(Y)|g(y)=f(x)\}$, and $p_1$ and $p_2$ are the projection digraph maps.
\end{definition}
\begin{proposition}\label{pullback}
For any commutative diagram
$$
\xymatrix{
   & G \ar[d]_{l_1} \ar[r]^{l_2}
                      & Y \ar[d]^{g}    \\
   & X \ar[r]^{f}     & Z,               }
$$ there exists a unique digraph map $p\colon  G\rightarrow X\times_{Z} Y$ such that $p_2\circ p = l_2$ and $p_1\circ p = l_1,$ that is, there is a commutative diagram of digraph maps:
$$
\xymatrix{
 G \ar@/_/[ddr]_{l_1} \ar@/^/[drr]^{l_2}
    \ar@{.>}[dr]|-{p}                   \\
   & X\times_{Z} Y \ar[d]^{p_1} \ar[r]_{p_2}\pushoutcorner
                      & Y \ar[d]_{g}    \\
   & X \ar[r]^{f}     & Z              . }
$$
\begin{proof}
Since $f\circ l_1 = g\circ l_2$, $(l_1(v),l_2(v))\in X\times_{Z}Y $ for all vertex $v\in G$. Now let us define a map $$p\colon  V(G)\rightarrow V(X\times_{Z}Y),\text{ }v\mapsto (l_1(v),l_2(v)).$$ If $(v, w) \in A(X\times_{Z}Y)$, then $p((v,w)) = \left(\left(l_1\left(v\right),l_2\left(v\right)\right),\left(l_1\left(w\right),l_2\left(w\right)\right)\right)$.
Since $l_1$ is a digraph map, $l_1(v)=l_1(w)$ or $l_1(v)\rightarrow l_1(w)$. Similarly for $l_2$. Then   $\left(l_1\left(v\right),l_2\left(v\right)\right)=\left(l_1\left(w\right),l_2\left(w\right)\right)$ or  $\left(l_1\left(v\right),l_2\left(v\right)\right)\rightarrow\left(l_1\left(w\right),l_2\left(w\right)\right)$. Hence $p$ is a digraph map. Next we will show $p$ is unique. Suppose there is another digraph map $p'$ such that $ p_1\circ p' = l_1$ and $p_2\circ p' = l_2$. Then $l_1(v) =  p_1( p'(v))$ and $l_2(v)=p_2(p'(v))$. Since $p'$ is a digraph map, it is determined the the map of vertices, therefor $p'=p$. That is to say, $p$ is unique.
\end{proof}
\end{proposition}
From the definition of pullback in the digraph category, the mapping path space $P_f$ for any based digraph map $f$ can be given by the following definition.
\begin{definition}Let $X$ be a digraph with base-point $x_0$ and $G$ be a digraph with base-point $\ast$, for any based digraph map $f\colon X\rightarrow G$, the \emph{mapping path digraph} $P_f$ is the pullback of $f$ and $e$, where $e$ is the evaluation map $e\colon  \overline{P}G\rightarrow G$. Explicitly, $P_f$ is the induced sub-digraph of $ X\times \overline{P}G$ with base-point $(x_0,\langle l_\ast\rangle)$ generated by $V(P_f)=\{(x,\langle\lambda\rangle)\}\in V(X)\times V(\overline{P}G^{\ast})|e(\langle\lambda\rangle)=f(x)\}$.
$$
\xymatrix{
  P_f \ar[d]_{f'} \ar[r]_{q} \pushoutcorner
                & \overline{P}G \ar[d]^{e}  \\
  X  \ar[r]^{f}
               & G             }
$$
\end{definition}
Following the idea from topological spaces, we verify a Puppe sequence of digraphs by using $P_f$.
\begin{proposition}
For any based digraph map $f\colon  X\rightarrow G$, the following sequence is exact at $\overline{\pi}_n(X)$ as set
$$\xymatrix@C=0.5cm{
  \overline{\pi}_n(P_f) \ar[r]^{f'_{n}} & \overline{\pi}_n(X) \ar[r]^{f_{n}} & \overline{\pi}_n(G),}$$
for $n\geq0$. When $n\geq 1$, it is an exact sequence of groups.
\begin{proof}
First of all, the maps $f'_{n}$ and $f_{n}$ are well-defined by Proposition \ref{induce} for $n\geq 0$. Since $f_n$ and $f'_n$ are homomorphisms if $n\geq1$, we only need to check it is an exact sequence. Now let us check $\im f'_{n}\subset \ker f_{n}$, that is, $f_{n} \circ f'_{n} =0$. Since $f\circ f'= e\circ q $, we obtain $f_{n} \circ f'_{n} =e_{n}\circ q_{n}, $ that is, there is a commutative diagram:
$$
\xymatrix{
   \overline{\pi}_n(P_f) \ar[d]_{f'_{n}} \ar[r]^{q_{n}}
                &  \pi_n(\overline{P}G) \ar[d]^{e_{n}}  \\
 \overline{\pi}_n(X)  \ar[r]^{f_{n}}
               & \overline{\pi}_n(G)           . }
$$
Since $\overline{P}G$ is weakly contractable by Proposition \ref{weak}, $f_{n} \circ f'_{n}=0$.

Next let us check $\ker f_{n}\subset \im f'_{n}.$ If $n\geq 1$, for any $[g]\in \overline{\pi}_n(X)$ such that $f \circ g\s_F \ast$, we will show that there exists a digraph map $\widetilde{f}\colon   (J_{M_i},\partial J_{M_i})^{\Box n}\rightarrow (P_f,(x_0,\langle l_\ast\rangle))$ such that $f'\circ \widetilde{f}\s_F g$ via a commutative diagram of digraph maps:
$$
\xymatrix{
 J_{M_i}^{\Box n} \ar[ddr]^{\overline{g}} \ar[d]_{q_l}  \ar@/^/[drr]^{\gamma}
    \ar@{.>}[dr]|-{\widetilde{f}}                   \\
J_{m_i}^{\Box n}\ar[dr]_{g}   &P_f \ar[d]^{f'} \ar[r]_{q}\pushoutcorner
                      &\overline{P}G \ar[d]^{e}    \\
   & X \ar[r]^{f}     & G               }
$$

Since $f \circ g\s_F \ast$, there is a finite sequence of digraph maps $\{h_i\}_{i=0}^{l}$ from $ \ast$ to $f \circ g$ such that $h_i\s_1 h_{i+1}$ or $h_i\s_{-1}h_{i+1}$, where $h_i\colon  (J_{m_j^i},\partial J_{m_j^i})^{\Box n}\rightarrow (G,\ast)$ for $0\leq i\leq  l-1$. Suppose $h_0\s_1 h_{1}$, then there exist subdivisions $\overline{h_0}$ and $\overline{h_{1}}$ of $h_0$ and $h_{1}$ respectively such that $\overline{h_0} \dr \overline{h_{1}}$ by $q_0$ and $q_{1}$. Next by Corollary \ref{subren}, there exist subdivisions $\overline{\overline{h_{1}}}$ and $\overline{h_2}$ of $\overline{h_{1}}$ and $h_2$ respectively by $q_1'$ and $q_2$. At the same time, we subdivide $\overline{h_0}$ by $q_1'$. Repeat this procedure iteratively for $0\leq i\leq l$ until we obtain a digraph map $$H\colon  (J_{M_i},\partial J_{M_i})^{\Box n}\Box (I_{l}, 0)\rightarrow (G,\ast)$$ such that $H|_{(J_{M_i},\partial J_{M_i})^{\Box n}\Box \{l\}} = h_l.$

 By Proposition \ref{duli}, there is a digraph map $\phi(H)\colon (J_{M_i},\partial J_{M_i})^{\Box n}\rightarrow ((G,\ast)^{(I_{l},0)},\ast^{I_{l}})$. Since $(G,\ast)^{(I_{l},0)}$ is a sub-digraph of $PG$ and there is a digraph map $p\colon  PG\rightarrow \overline{P}G$, we obtain a digraph map
$$ \gamma = p\circ i \circ \phi(H)\colon (J_{M_i},\partial J_{M_i})^{\Box n}\rightarrow (\overline{P}G,\langle l_\ast\rangle).$$

It is clear that $e\circ\gamma(a) = \overline{h_l}(a)$ for any $a\in J_{M_i}^{\Box n}$, where $\overline{h_l}$ is the subdivision of $f \circ g$ by $q_l$, that is, $\overline{h_l} =f \circ g\circ q_l$. Or we can write $\overline{h_l} =f \circ \overline{g}$, where $\overline{g}$ is the subdivision of $g$ by $q_l$. Hence we have $e\circ\gamma= f \circ\overline{g}$. By Proposition \ref{pullback}, there exists a unique digraph map $\widetilde{f}\colon  J_{M_i}^{\Box n}\rightarrow P_f$ such that $\overline{g}= f'\circ \widetilde{f}$. Moreover $\widetilde{f} $ is a relative digraph map  $(J_{M_i},\partial J_{M_{i}})^{\Box n}\rightarrow (P_f,(\ast,\langle l_\ast\rangle)) .$ So $f'_{n}( [\widetilde{f}])=[\overline{g}]=[g].$ Thus $\ker f_{n}\subset \im f'_{n}.$

If $n=0,$  for any $g\in \operatorname{Hom}(0,X)$ such that $ f \circ g\s_F \ast$, there exists a line digraph $I_l$ and a digraph map $F\colon I_l\rightarrow G$ such that $F(0) = \ast(0)=\ast$ and $F(l) = f\circ g(0)$. Then we can construct two digraph maps $H\colon  0\rightarrow PG$ such that $H(0)=F$ and $\gamma := p\circ H\colon  0\rightarrow \overline{P}G$. Since $e\circ \gamma = f \circ g$, by the universal property of a pullback, there exists a unique digraph map $\widetilde{f}\colon  0 \rightarrow P_f$ such that $f'\circ \widetilde{f} = g$, that is, there is a commutative diagram:
$$
\xymatrix{
 0 \ar@/_/[ddr]^{g}  \ar@/^/[drr]^{\gamma}
    \ar@{.>}[dr]|-{\widetilde{f}}                   \\
  &P_f \ar[d]^{f'} \ar[r]_{q}\pushoutcorner
                      &\overline{P}G \ar[d]^{e}    \\
   & X \ar[r]^{f}     & G               }
$$ Hence $f'_{0}([\widetilde{f}]) = [g]$. Thus $\ker f_{0}\subset \im f'_{0}.$
In conclusion, $$\xymatrix@C=0.5cm{
  \overline{\pi}_n(P_f) \ar[r]^{f'_{n}} & \overline{\pi}_n(X) \ar[r]^{f_{n}} & \overline{\pi}_n(G) }$$
 is an exact sequence for $n\geq 0.$
\end{proof}
\end{proposition}
Based on this exact sequence, we obtain the following commutative diagram by using the pullback of digraphs iteratively:
$$
\xymatrix{
  P_{f^{(3)}} \ar[d]_{f^{(4)}} \ar[r] \pushoutcorner & \overline{P}(P_{f'}) \ar[d] &  \\
  P_{f^{(2)}} \ar[d] \ar[r]^{f^{(3)}}\pushoutcorner & P_{f'} \ar[d]_{f^{(2)}} \ar[r]\pushoutcorner & \overline{P}X \ar[d] \\
  \overline{P}(P_{f}) \ar[r] & P_{f} \ar[d] \ar[r]^{f'} \pushoutcorner& X \ar[d]^{f} \\
  & \overline{P}G \ar[r] & G.   }
  $$
Further, it induces a long exact sequence of homotopy groups.
\begin{corollary}\label{les}
For any based digraph map $f\colon  X\rightarrow G$, there is a long exact sequence
$$\xymatrix@C=0.5cm{
 \ar[r]& \cdots \ar[r]& \overline{\pi}_n (P_{f^{(2)}}) \ar[r]^{f^{(3)}_{n}}& \overline{\pi}_n(P_{f'}) \ar[r]^{f^{(2)}_{n}}&
  \overline{\pi}_n(P_f) \ar[r]^{f'_{n}}&
 \overline{\pi}_n (X) \ar[r]^{f_{n}} & \overline{\pi}_n(G) },$$
of based sets for $n\geq0$. If $n\geq 1$, it is a long exact sequence of groups.\qed
\end{corollary}
\begin{proposition}\label{Pf}
There is an isomorphism $\iota\colon  P_{f'}\rightarrow H_{f}$, where $H_f$ is the induced sub-digraph of $\overline{P}G\times \overline{P}X$ generated by $V(H_{f}) = \{(\langle\gamma\rangle,\langle\eta\rangle)|f(e_{X}(\langle\eta\rangle)) =e_{G}(\langle\gamma\rangle) \}$.
\begin{proof}
As we mentioned, $P_{f'}$ is a digraph with vertex set $$V(P_{f'})=\{(x, \langle\gamma\rangle,\langle\eta\rangle)\in X\times \overline{P}G\times \overline{P}X| f(x)=e_{G}(\langle\gamma\rangle), e_{X}(\langle\eta\rangle)=x\}.$$
Define a map $$\iota\colon P_{f'}\rightarrow H_{f}$$ by $$(x, \langle\gamma\rangle,\langle\eta\rangle)\mapsto (\langle\gamma\rangle,\langle\eta\rangle).$$

Clearly, $\iota$ is an isomorphism of digraphs.
\end{proof}
\end{proposition}
According to the definition of $P_{f'}$, there is a based digraph map $$j\colon  (\overline{L}G,\langle l_\ast\rangle)\rightarrow \left(P_{f'},\left(x_0,\langle l_\ast\rangle,\langle l_{x_0}\rangle\right)\right),\text{ }\langle\gamma\rangle\mapsto (x_0,\langle\gamma\rangle,\langle l_{x_0}\rangle),$$ where $l_{x_0}\colon (J_1,\partial J_1)\rightarrow x_0$ and $l_{\ast}\colon  (J_1,\partial J_1)\rightarrow \ast.$ From $j$, we obtain the following theorem.
\begin{theorem}\label{nL}The digraph map $j$ induces an isomorphism $$j_{n}\colon \overline{\pi}_n(\overline{L}G) \rightarrow \overline{\pi}_n(P_{f'})$$ for $n\geq0$.
\begin{proof}
By Proposition \ref{induce}, a digraph map induces a morphism between homotopy groups. To prove the theorem, it is sufficient to show that there exists a digraph map $$q\colon P_{f'}\rightarrow \overline{L}G,\text{ }(x,\langle\gamma\rangle,\langle\eta\rangle)\mapsto \langle \gamma \vee(f\circ \eta^{-1})\rangle$$ such that $j_n\circ q_n= id_{\overline{\pi}_{n}(P_{f'})}$ and $q_n\circ j_n = id_{\overline{\pi}_{n}(\overline{L}G)},$ where $\eta^{-1}$ is the inverse digraph map of $\eta\colon  (J_{m},0)\rightarrow (X,x_0)$.

Let $\overline{\gamma}\colon  (J_M,0)\rightarrow (G,\ast)$ be a subdivision of $\gamma$ by $q_1$. Given any subdivision $\overline{\eta}\colon (J_N,0)\rightarrow (X,x_0)$ of $\eta$ by $q_2$, one can easily show that the inverse map $\overline{\eta}^{-1}$ of $\overline{\eta}$ is the subdivision of $\eta^{-1}$ by $q_2^{-1}$, so $\langle\overline{\eta}^{-1}\rangle= \langle\eta^{-1}\rangle$. Then $\overline{\gamma}\vee(f\circ\overline{\eta}^{-1})$ is the subdivision of $\gamma\vee(f\circ \eta^{-1})$ by $q_1\vee q_2^{-1}$, where $q_1\vee q_2^{-1}\colon I_{M+N}\rightarrow I_{m+n}$ is the concatenation of $q_1$ and $q_2^{-1}$. Hence the multiplication and the inverse of maps is independent of the subdivision. For any vertex $\langle\gamma\rangle$ in $\overline{P}G$, there exists a unique minimal path $\widehat{\gamma}$ in $\langle\gamma\rangle$ by deleting adjacent vertices that are the same, and any element in $\langle\gamma\rangle$ is a subdivision of $\widehat{\gamma}$. Thus $q$ is well-defined on the vertex set.

Assume  $(x,\langle\gamma\rangle,\langle\eta\rangle)\rightarrow(x',\langle\gamma'\rangle,\langle\eta'\rangle)$ in $P_{f^{'}}$, even if $\langle\gamma\rangle\rightarrow \langle\gamma'\rangle$ and $\langle\eta\rangle\rightarrow \langle\eta'\rangle$, there exist subdivisions $\overline{\gamma}$ of $\gamma$, $\overline{\gamma'}$ of $\gamma'$, $\overline{\eta}$ of $\eta$ and $\overline{\eta'}$ of $\eta'$ by $p_1$, $p'_1$, $p_2$ and $p_2'$ respectively such that $\overline{\gamma}\dr \overline{\gamma'}$ and $\overline{\eta}\dr \overline{\eta'}$. It follows that $\overline{\gamma}\vee (f\circ\overline{\eta}^{-1})\dr \overline{\gamma'}\vee (f\circ(\overline{\eta'})^{-1})$. Hence $\langle\gamma\vee (f\circ \eta^{-1})\rangle\rightarrow\langle\gamma'\vee (f\circ(\eta')^{-1})\rangle$. Thus $q$ is a digraph map. Further $q$ preserves the base-point. By Proposition \ref{induce}, $q$ induces a homomorphism
 $$ q_n\colon \overline{\pi}_n(P_{f'})\rightarrow \overline{\pi}_n(\overline{L}G)$$
for $n\geq 0.$

Consider $n\geq 1.$ For any $h\colon (J_{m_i},\partial J_{m_i})^{\Box n}\rightarrow (\overline{L}G, \langle l_\ast \rangle)$ and  $a\in (J_{m_i},\partial J_{m_i})^{\Box n}$, $$q\circ j\circ h(a) = q\circ j(\langle\gamma_a\rangle) = q((x_0,\langle\gamma_a\rangle,\langle l_{x_0}\rangle))= \langle\gamma_a\vee (f\circ l_{x_0}^{-1})\rangle =\langle\gamma_a\rangle= h(a),$$
where $h(a) = \langle\gamma_a\rangle$. So $q_n\circ j_n([h])=[q\circ j\circ h]=[h]$. Hence $q_n\circ j_n = id_{\overline{\pi}_{n}(\overline{L}G)}.$ On the other hand, we claim $j_n\circ q_n = id_{\overline{\pi}_{n}(P_{f^{'}})}$, that is, $j\circ q\circ h\s_F h$ for any $[h]\in \overline{\pi}_n(P_{f'})$. In our case, we may assume that $h$ come from a map $$h\colon (J_{m_i},\partial J_{m_{i}})^{\Box n}\rightarrow (P_{f'},(x_0,\langle l_{\ast}\rangle,\langle l_{x_0}\rangle))$$ with $h(a) = (x_a,\langle\gamma_a\rangle,\langle\eta_a\rangle)$. Then $j\circ q\circ h(a) = (x_0, \langle\gamma_a\vee (f\circ \eta_a^{-1})\rangle,\langle l_{x_0}\rangle)$ for any $a\in J_{m_i}^{\Box n}$. If we focus on the $\langle\eta_a\rangle$ part in $P_{f'}$, there is a digraph map $\overline{h}\colon (J_{m_i},\partial J_{m_i})^{\Box n}\Box  $ $(J_M,0)\rightarrow(X,x_0)$ such that $\overline{h}(a,-) =\overline{\eta_a} $ for all  $a \in J_{m_i}^{\Box n}$ by Remark \ref{plattice}, where $\overline{\eta_a} $ is a subdivision of $\eta_a$. Then we choose $\overline{\eta_a}\colon (J_{M},0)\rightarrow (X,x_0)$ as the representative element of $\langle\eta_a\rangle$. Now we can construct a sequence of digraph maps $\{g_j \}_{j=0}^{M}$, where $$g_j\colon (J_{m_i},\partial J_{m_{i}})^{\Box n}\rightarrow (P_{f'},(x_0,\langle l_{\ast}\rangle,\langle l_{x_0}\rangle))$$ is defined by $$i\mapsto (\overline{\eta_i}|_{J_{M-j}}(M-j),\langle\gamma_a\vee (f\circ (\overline{\eta_{a}}^{-1}|_{J_j} )) \rangle,\langle\overline{\eta_i}|_{J_{M-j}}\rangle).$$
Clearly, $g_M= j\circ q\circ h$ and $g_0 = h$. If $M-j-1\rightarrow M-j$, then $\overline{\eta_i}|_{J_{M-j}}\s_{-1} \overline{\eta_i}|_{J_{M-j-1}}$ and $\overline{\eta_{i}}^{-1}|_{J_j} \s_{-1}\overline{\eta_{i}}^{-1}|_{J_{j+1}} $ for any $i$. Hence $g_{j}\dl g_{j+1}$. Similarly, if $M-j-1\leftarrow M-j$, then $g_{j}\dr g_{j+1}$. Thus $j\circ q\circ h\s_F h$. Since $h$ is arbitrary, we obtain $j_n\circ q_n = id$. Hence $j_n$ is an isomorphism for $n\geq 1$.

For $n=0,$ let $h\colon 0\rightarrow \overline{L}G$ be any digraph map such that $h(0)=\langle\gamma\rangle$. Since $q\circ j\circ h(0) = h(0)=\langle\gamma\rangle$, $q_0\circ j_0=id_{\overline{\pi}_0(\overline{L}G)}$. Now let us check $j_0\circ q_0 =id_{\overline{\pi}_0(P_{f'})}.$ Given any digraph map $h\colon 0\rightarrow P_{f'}$ such that $h(0)=(x,\langle\gamma\rangle,\langle\eta\rangle)$, we have $j\circ q\circ h(0)= (x_0,\langle\gamma\vee (f\circ\eta^{-1})\rangle,\langle l_{x_0}\rangle)$. It remains to prove that there is a digraph map $F\colon I_m\rightarrow P_{f'}$ such that $F(0)=(x,\langle\gamma\rangle,\langle\eta\rangle)$ and $F(m)=(x_0,\langle\gamma\vee (f\circ\eta^{-1})\rangle,\langle l_{x_0}\rangle).$

We construct a map $$F\colon I_m\rightarrow P_{f'}, \text{ }j\mapsto (\eta|_{J_{m-j}}(m-j),\langle\gamma\vee (f\circ (\eta^{-1}|_{J_j} )) \rangle,\langle\eta|_{J_{m-j}}\rangle)$$
such that $F(m)= j\circ q\circ h(0)$ and $F(0) = h(0)$. Since $\eta|_{J_{m-j}}\s_1\eta|_{J_{m-j+1}}$ or $\eta|_{J_{m-j}}\s_{-1}\eta|_{J_{m-j+1}}$, then $f\circ (\eta^{-1}|_{J_j})\s_1 f\circ (\eta^{-1}|_{J_{j-1}})$ or $f\circ (\eta^{-1}|_{J_j})\s_{-1} f\circ (\eta^{-1}|_{J_{j-1}})$ respectively. By the Cartesian product of digraphs, $F$ is a digraph map. Hence $j_0\circ q_0([h])=[j\circ q\circ h]=[h]=id([h])$ for any $[h]$. Because $h$ is arbitrary, we obtain $j_0\circ q_0 =id_{\overline{\pi}_0(P_{f'})}.$
\end{proof}
\end{theorem}

By Theorem \ref{nL}, if we replace $G$ by $X$ and replace $P_{f}$ by $P_{f^{(2)}}$, then  we get $\overline{\pi}_n(\overline{L}X)\stackrel{j_n}{\approx} \overline{\pi}_n(P_{f^{(2)}}).$

Since isomorphisms cannot affect the exactness of the exact sequence, we obtain a long exact sequence \begin{align}\label{replace}
\quad
\quad \quad
\quad
\quad
\xymatrix@C=0.5cm{
    \cdots& \ar[r] &\overline{\pi}_{n+1}(X) \ar@{.>}[r]^{\overline{f^{(3)}_{n}}} \ar[d]^{\Phi^{n+1}_{\ast}}&  \overline{\pi}_{n+1}(G) &
    &
    & \\
  \cdots&\ar[r]& \overline{\pi}_{n}(\overline{L}X) \ar[d]^{j_n}& \overline{\pi}_{n}(\overline{L}G)\ar[u]^{\Psi^{n+1}_{\ast}} &
    &
    & \\
\cdots&  \ar[r] &\overline{\pi}_{n}(P_{f^{(2)}}) \ar[r]^{f^{(3)}_{n}}&  \overline{\pi}_{n}(P_{f'})\ar@<.5ex>[u]^{q_n}  \ar[r]^{f^{(2)}_{n}}&
    \overline{\pi}_n(P_f)\ar[r]^{f'_{n}}&
  \overline{\pi}_n(X) \ar[r]^{f_{n}} & \overline{\pi}_n(G)}
  \end{align}
 by using isomorphisms iteratively for any $n\geq0$, where $\overline{f^{(3)}_{n}}= \Psi^{n+1}_{\ast}\circ q_n \circ f^{(3)}_{n}\circ j_n\circ \Phi^{n+1}_{\ast}$. Then we write $\overline{f_{n}^{(3)}}$ as $\Omega f_n$.

We finally get our Puppe sequence of digraphs by considering the morphisms in diagram \ref{replace}.
\begin{theorem}[\bf{Puppe Sequence}]\label{puppe}
For any based digraph map $f\colon  X\rightarrow G$, there is a long exact sequence
$$\xymatrix@R=0.4cm{
 \cdots\ar[r]&\overline{\pi}_{n+2}(X) \ar[r]^{ f_{n+2}}&\overline{\pi}_{n+2}(G)\ar[r]^{\Omega\partial_{n+1}}&    \overline{\pi}_{n+1}(P_f) \ar[r]^{\Omega f^{'}_{n}}& \overline{\pi}_{n+1}(X) \\ \ar[r]^-{\Omega f_{n}}&  \overline{\pi}_{n+1}(G) \ar[r]^{\partial_{n+1}}&
    \overline{\pi}_n(P_f)\ar[r]^{f^{'}_{n}}&
  \overline{\pi}_n(X) \ar[r]^{f_{n}} & \overline{\pi}_n(G)}$$
of based sets for any $n\geq0 $. If $n\geq 1$, it is a long exact sequence of groups.
   \begin{proof}
First, let us show that there is a commutative diagram,
$$\xymatrix{
  P_{f^{(2)}}  \ar[r]^{f^{(3)}}
                & P_{f'} \ar[d]^{q_{\overline{L}G}}  \\
  \overline{L}X \ar[u]^{j_{\overline{L}X}} \ar@{.>}[r]^{\overline{L}f\circ \nu}
                & \overline{L}G             ,}$$
where $$\nu\colon  \overline{L}X\rightarrow \overline{L}X, \quad \langle\varepsilon\rangle\mapsto \langle\varepsilon^{-1}\rangle.$$
For any $\langle\varepsilon\rangle\in V(\overline{L}X)$,\begin{align*}
q_{\overline{L}G}\circ f^{(3)}\circ j_{\overline{L}X} (\langle\varepsilon\rangle) &= q_{\overline{L}G}\circ f^{(3)}(\ast_{P_f},\langle\varepsilon\rangle,\langle\ast\rangle_{\overline{P}(P_f)})\\
 &=q_{\overline{L}G}((\ast_{P_f},\langle\varepsilon\rangle))\\
&= \langle l_{\ast}\vee (f\circ \varepsilon^{-1})\rangle\\
&= \langle f\circ \varepsilon^{-1}\rangle.
\end{align*}
On the other hand, it is easy to see that $\nu$ is a digraph map and $\overline{L}f\circ \nu(\langle\varepsilon\rangle) = \langle f\circ \varepsilon^{-1}\rangle = q_{\overline{L}G}\circ f^{(3)}\circ j_{\overline{L}X} (\langle\varepsilon\rangle)$ for any $\langle\varepsilon\rangle\in V(\overline{L}X)$. So $q_{\overline{L}G}\circ f^{(3)}\circ j_{\overline{L}X}=\overline{L}f\circ \nu.$

We turn to consider the morphisms at the level of homotopy groups because $\overline{L}G$ is weakly homotopy equivalent to $P_{f'}$ instead of homotopy equivalent. Assume $n\geq 1$. Consider the diagram
$$
\xymatrix{
  \overline{\pi}_{n}(P_{f^{(2)}}) \ar[r]^{f^{(3)}_{n}} & \overline{\pi}_{n}(P_{f^{'}}) \ar@<.5ex>[d]^{q_n}\\
  \overline{\pi}_{n}(\overline{L}X) \ar@<.5ex>[u]^{j_n}\ar[r]^{\overline{L}f_{n}\circ \nu_{n}} & \overline{\pi}_{n}(\overline{L}G) \ar[d]^{\Psi^{n+1}_{\ast}} \\
  \overline{\pi}_{n+1}(X) \ar[u]^{\Phi^{n+1}_{\ast}}\ar@{.>}[r]^{\Omega f_n} & \overline{\pi}_{n+1}(G)
   }
$$
The top square commute by the previous paragraph. We will define the map $\Omega f_n$
to make the lower square commute.

For any $\gamma\colon  (J_{m_i}, \partial J_{m_i})^{\Box (n+1)}\rightarrow (X,x_0)$,
\begin{align*}
\Psi^{n+1}_{\ast}\circ q_n\circ f^{(3)}_{n}\circ j_n\circ \Phi^{n+1}_{\ast}([\gamma])&= \Psi^{n+1}_{\ast}\circ \overline{L}f_{n} \circ \nu_{n}\circ \Phi^{n+1}_{\ast}([\gamma])\\
&=[\overline{L}f\circ \nu\circ \widetilde{\gamma}]\\&= [f\circ \gamma^{-1}],
\end{align*}
where $\widetilde{\gamma} = \phi_{m_{n+1}}(\gamma).$
Denoted $\Psi^{n+1}_{\ast}\circ q_n\circ f^{(3)}_{n}\circ j_n\circ \Phi^{n+1}_{\ast}([\gamma])$ by $\Omega f_{n}([\gamma])$. Similarly, $\Omega f^{(3)}_{n}([\gamma])= [f_n^{(3)}\circ\gamma^{-1}] $, where $\gamma\colon  (J_{m_i},\partial J_{m_i})^{\Box n}\rightarrow (P_{f^{(2)}},\ast)$.

Further, let us compute $\Omega^{2}f_{n}:=\Psi^{n+2}_{\ast}\circ q_{n+1}\circ \Psi^{n+1}_{\ast}\circ q_n\circ f^{(6)}_{n}\circ j_n\circ \Phi^{n+1}_{\ast}\circ j_{n+1}\circ \Phi^{n+2}_{\ast} $. Consider the diagram
$$
\xymatrix{   \overline{\pi}_{n}(P_{f^{(5)}})  \ar[r]^{f^{(6)}_{n}} & \overline{\pi}_{n}(P_{f^{(4)}}) \ar@<.5ex>[d]^{q_n}\\
 \overline{\pi}_{n}(\overline{L}P_{f^{(2)}}) \ar@<.5ex>[u]^{j_n}\ar[r]^{\overline{L}f^{(3)}_{n}\circ \nu_{n}} & \overline{\pi}_{n}(\overline{L}P_{f^{'}}) \ar[d]^{\Psi^{n+1}_{\ast}} \\
  \overline{\pi}_{n+1}(P_{f^{(2)}}) \ar[u]^{\Phi^{n+1}_{\ast}} \ar[r]^{\Omega f^{(3)}_{n}} & \overline{\pi}_{n+1}(P_{f^{'}}) \ar@<.5ex>[d]^{q_{n+1}}\\
  \overline{\pi}_{n+1}(\overline{L}X) \ar@<.5ex>[u]^{j_{n+1}}& \overline{\pi}_{n+1}(\overline{L}G) \ar[d]^{\Psi^{n+2}_{\ast}} \\
  \overline{\pi}_{n+2}(X) \ar[u]^{\Phi^{n+2}_{\ast}} \ar@{.>}[r]^{\Omega^{2}f_{n}}& \overline{\pi}_{n+2}(G)
   }
 $$
Assume $\eta\colon (J_{m_i},\partial J_{m_i})^{\Box(n+2)}\rightarrow (X,x_0)$ is defined by $\eta(a,i,j)= x_{a,i,j}$. Since $ \Omega f^{(3)}_{n}= \Psi^{n+1}_{\ast}\circ q_n\circ f^{(6)}_{n}\circ j_n\circ \Phi^{n+1}_{\ast},$  \begin{align*}
\Omega^{2}f_{n}([\eta]) = \Psi^{n+2}_{\ast}\circ q_{n+1}\circ\Omega f^{(3)}_{n}\circ j_{n+1}\circ \Phi^{n+2}_{\ast}
= [f \circ \widehat{\eta}]
=f_{n+2}([\widehat{\eta}]),
\end{align*}
where $\widehat{\eta}(a,i,j) = x_{a,m_{n+1}-i,m_{n+2}-j}$, $a\in J_{m_{i}}^{\Box n}$, $i\in J_{m_{n+1}}$ and $j\in J_{m_{n+2}}$.

By the description of the inverse element of $\eta$ in Theorem \ref{gp}, $\eta_{n+2}^{-1}(a,i,j) = x_{a,i,m_{n+2}-j}$ and $ (\eta_{n+2}^{-1})_{n+1}^{-1}(a,i,j) = x_{a,m_{n+1}-i,m_{n+2}-j}$, so $\widehat{\eta} = (\eta_{n+2}^{-1})_{n+1}^{-1}$, therefore $[\widehat{\eta}] = [(\eta_{n+2}^{-1})]^{-1} = ([\eta]^{-1})^{-1} = [\eta].$ Hence $\Omega^{2}f_{n}([\eta]) = f_{n+2}([\widehat{\eta}]) = f_{n+2}([\eta]) .$ Thus $\Omega^{2}f_{n}$ $ = f_{n+2}$.

Now consider $n=0.$ By definition of the fundamental groups of digraphs, there is an isomorphism $$\Theta_0\colon \overline{\pi}_1(X)\rightarrow \overline{\pi}_0(LX),\quad [\eta]\mapsto [\widetilde{\eta}],$$ where $\widetilde{\eta}\colon  0\rightarrow LX$ such that $\widetilde{\eta}(0)=\eta\colon (J_m,\partial J_m)\rightarrow (X,x_0)$. The based digraph map $p\colon  (LX,l_{\ast})\rightarrow (\overline{L}X, \langle l_{\ast}\rangle)$ induces an isomorphism $p_0$ by Proposition \ref{0i}, so we define $\Phi^1_{\ast} = p_0\circ \Theta_0\colon  \overline{\pi}_1(X)\rightarrow \overline{\pi}_0(\overline{L}X).$ Consider the diagram
   $$
\xymatrix{
  \overline{\pi}_{0}(P_{f^{(2)}})  \ar[r]^{f^{(3)}_{0}} & \overline{\pi}_{0}(P_{f^{'}}) \ar@<.5ex>[d]^{q_0}\\
  \overline{\pi}_{0}(\overline{L}X) \ar@<.5ex>[u]^{j_0}\ar[r]^{\overline{L}f_{0}\circ \nu_{0}} & \overline{\pi}_{0}(\overline{L}G) \ar[d]^{\Psi^1_{\ast} } \\
  \overline{\pi}_{1}(X) \ar[u]^{\Phi^1_{\ast} }\ar@{.>}[r]^{\Omega f_{0}} & \overline{\pi}_{1}(G).
   }
$$
By computing, for any $ [\eta]\in \overline{\pi}_1(X)$, $\Omega f_{0}([\eta]) = [f\circ \eta^{-1}].$  Similarly to the $n\geq 1$ case, $\Omega^{2} f_{0}= f_{2}.$

Now we define the homomorphisms $\partial_{n+1}$ for $n\geq 1$ by the composite
$$
\xymatrix{
  \overline{\pi}_{n}(P_{f^{'}}) \ar[r]^{f^{(2)}_{n}} & \overline{\pi}_{n}(P_{f}) \\
  \overline{\pi}_{n}(\overline{L}G) \ar@<.5ex>[u]^{j_n} &    \\
  \overline{\pi}_{n+1}(G) \ar[u]^{\Phi^{n+1}_{\ast}}\ar@{.>}[uur]_{\partial_{n+1}}. &
}
$$
Assume $\gamma\colon  (J_{m_i} \partial J_{m_i})^{\Box (n+1)}\rightarrow (G,\ast)$. Then $\widetilde{\gamma}=\phi_{m_{n+1}}(\gamma)\colon (J_{m_i}, \partial J_{m_i})^{\Box n}\rightarrow (\overline{L}G,\langle l_\ast\rangle)$, and $$\partial_{n+1}([\gamma]):=f^{(2)}_{n}\circ j_n\circ\Phi^{n+1}_{\ast}([\gamma])=[f^{(2)}\circ j\circ \widetilde{\gamma}].$$ Consider the diagram
 $$
\xymatrix{
\overline{\pi}_{n}(P_{f^{(4)}})\ar[r]^{f^{(5)}_{n}} & \overline{\pi}_{n}(P_{f^{(3)}}) \ar@<.5ex>[d]^{q_n}\\
  \overline{\pi}_{n}(\overline{L}P_{f'})\ar[u]^{j_{n}} & \overline{\pi}_{n}(\overline{L}P_{f}) \ar[d]^{\Psi^{n+1}_{\ast}} \\
  \overline{\pi}_{n+1}(P_{f^{'}}) \ar[u]^{\Phi^{n+1}_{\ast}} \ar[r]^{\Omega f^{(2)}_{n}} & \overline{\pi}_{n+1}(P_{f}) \\
  \overline{\pi}_{n+1}(\overline{L}G) \ar@<.5ex>[u]^{j_{n+1}} &    \\
  \overline{\pi}_{n+2}(G). \ar[u]^{\Phi^{n+2}_{\ast}}\ar@{.>}[uur]_{\Omega\partial_{n+1}} &.
   }
$$
For any $ [\eta]\in \overline{\pi}_{n+1}(P_{f^{'}})$, we have $\Omega f^{(2)}_{n}([\eta]) = [f^{(2)}\circ \eta^{-1}].$ Therefore for any \\ $\gamma\colon  (J_{m_i}, \partial J_{m_i})^{\Box (n+2)}\rightarrow (G,\ast)$, we have $\widetilde{\gamma}:=\phi_{m_{n+2}}(\gamma)\colon (J_{m_i}, \partial J_{m_i})^{\Box (n+1)}\rightarrow (\overline{L}G,\langle l_\ast\rangle)$, $\Omega\partial_{n+1}([\gamma]):=\Omega f^{(2)}_{n}\circ j_{n+1}\circ\Phi^{n+2}_{\ast}([\gamma])=[f^{(2)}\circ (j\circ \widetilde{\gamma})^{-1}].$ Furthermore, as $\Omega^{2} f_{n}= f_{n+2},$ we obtain $\Omega^{2}\partial_{n+1}=\partial_{n+3}$.

Similarly, we compute the $\partial_1 = f^{(2)}_{0} \circ j_0\circ \Phi^{1}_{\ast}.$ For any $[\eta]\in \overline{\pi}_1(G)$, $\partial_1([\eta])=f^{(2)}_0 \circ j_0\circ \Phi^1_{\ast}([\eta]) = [f^{(2)}\circ j\circ \widetilde{\eta}]$, where $\widetilde{\eta}\colon  0\rightarrow \overline{L}G$ is the map such that $\widetilde{\eta}(0)=\eta\colon (J_m,\partial J_m) $ $ \rightarrow (G,\ast)$. Also, $\Omega\partial_1\colon \overline{\pi}_{2}(G)\rightarrow \overline{\pi}_1(P_f)$ is $\Omega\partial_1([\gamma])=\Omega f^{(2)}_0\circ j_1\circ \Phi^2_{\ast}([\gamma])$ $ = [f^{(2)}\circ (j\circ \widetilde{\gamma})^{-1}],$ where $\widetilde{\gamma} = \phi_{J_{m_2}}(\gamma).$  Moreover, as $\Omega^{2} f_{0}^{(2)}= f_{2}^{(2)},$ we obtain  $\Omega^{2}\partial_{1}=\partial_{3}$.

 Iteratively using the relation $\Omega^2 f_n^{(i)} = f_{n+2}^{(i)}$ for any $i\geq 1$, this theorem is proved.
   \end{proof}
   \end{theorem}
Here it should be pointed out that our Puppe sequence holds only for homotopy groups, instead of any set of homotopy class $[X,G]$ for a digraph $X$ since  $\overline{P}G$ is weakly contractible, which is different from the case in topological spaces. Nevertheless, we still have the following property as in classical homotopy theory.
\begin{proposition}
Any commutative diagram of based digraph maps $$\xymatrix{
  X \ar[d]_{u} \ar[r]^{f}
                & G \ar[d]^{v}  \\
  Y  \ar[r]^{g}
                & H             } $$induces a commutative diagram of exact sequences
$$\xymatrix{
\cdots\ar[r]&\overline{\pi}_{n+1}(X)\ar[r]^{\Omega f_n}\ar[d]_{u_{n+1}}&\overline{\pi}_{n+1}(G)\ar[d]^{v_{n+1}}\ar[r]^{\partial_n}&\overline{\pi}_n(P_f)\ar[r]^{f^{'}_{n}}\ar[d]^{w_{n}} & \overline{\pi}_n(X) \ar[d]_{u_{n}} \ar[r]^{f_{n}}
                & \overline{\pi}_n(G) \ar[d]^{v_{n}}  \\
\cdots\ar[r]&\overline{\pi}_{n+1}(Y)\ar[r]^{\Omega g_n}&\overline{\pi}_{n+1}(H)\ar[r]^{\partial_n}&\overline{\pi}_n(P_g)\ar[r]^{g^{'}_{n}} & \overline{\pi}_n(Y)  \ar[r]^{g_{n}}
                & \overline{\pi}_n(H)             }
$$ for $n\geq 0$.
\begin{proof}
Any based digraph map $v\colon  G\rightarrow H$ induces a digraph map $$\overline{v}\colon  \overline{P}G\rightarrow \overline{P}H, \quad \langle\gamma\rangle\mapsto \langle v\circ \gamma\rangle$$ such that $ e_{\overline{P}H}\circ \overline{v}=v\circ e_{\overline{P}G}.$ Therefore there is a digraph map $$w = u\times \overline{v}\colon  P_{f}\rightarrow P_{g} ,\quad (x,\langle\gamma\rangle)\mapsto (u(x), \langle v\circ\gamma\rangle)$$ such that
$g^{'}\circ w=u\circ f^{'}.$ Hence there is a commutative diagram
$$\xymatrix{
\overline{\pi}_n(P_f)\ar[r]^{f^{'}_{n}}\ar[d]_{w_{n}} & \overline{\pi}_n(X) \ar[d]_{u_{n}} \ar[r]^{f_{n}}
                & \overline{\pi}_n(G) \ar[d]^{v_{n}}  \\
\overline{\pi}_n(P_g)\ar[r]^{g^{'}_{n}} & \overline{\pi}_n(Y)  \ar[r]^{g_{n}}
                & \overline{\pi}_n(H)            . }$$
Iteratively, we obtain the following commutative diagram
$$\xymatrix{
\cdots\ar[r]&\overline{\pi}_n(P_{f^{'}})\ar[d]_{(w\times \overline{u})_n}\ar[r]^{f_n^{(2)}}&\overline{\pi}_n(P_f)\ar[r]^{f^{'}_{n}}\ar[d]^{w_{n}} & \overline{\pi}_n(X) \ar[d]_{u_{n}} \ar[r]^{f_{n}}
                & \overline{\pi}_n(G) \ar[d]^{v_{n}}  \\
\cdots\ar[r]&\overline{\pi}_n(P_{g^{'}})\ar[r]^{g_n^{(2)}}&\overline{\pi}_n(P_g)\ar[r]^{g^{'}_{n}} & \overline{\pi}_n(Y)  \ar[r]^{g_{n}}
                & \overline{\pi}_n(H)             }$$
for $n\geq 0$. There is also a commutative diagram
$$
\xymatrix{
  \overline{\pi}_{n}(P_{f^{'}})  \ar[r]^{(w\times \overline{u})_n} & \overline{\pi}_{n}(P_{g^{'}}) \ar@<.5ex>[d]^{q_n}\\
  \overline{\pi}_{n}(\overline{L}G) \ar@<.5ex>[u]^{j_n}& \overline{\pi}_{n}(\overline{L}H)\ar[d]^{\Psi^{n+1}_{\ast}} \\
  \overline{\pi}_{n+1}(G) \ar[u]^{\Phi^{n+1}_{\ast}}\ar[r]^{v_{n+1}} & \overline{\pi}_{n+1}(H),
   }
$$
since $$\Psi^{n+1}_{\ast}\circ q_n\circ (w\times \overline{u})_n \circ j_n\circ \Phi^{n+1}_{\ast}([\eta]) = [ v\circ\eta]= v_{n+1}([\eta]).$$
Iteratively, we obtain the following commutative diagram
$$\xymatrix{
\cdots\ar[r]&\overline{\pi}_{n+1}(X)\ar[d]_{u_{n+1}}\ar[r]^{\Omega f_n}&\overline{\pi}_{n+1}(G)\ar[d]^{v_{n+1}}\ar[r]^{\partial_n}&\overline{\pi}_n(P_f)\ar[r]^{f^{'}_{n}}\ar[d]^{w_{n}} & \overline{\pi}_n(X) \ar[d]_{u_{n}} \ar[r]^{f_{n}}
                & \overline{\pi}_n(G) \ar[d]^{v_{n}}  \\
\cdots\ar[r]&\overline{\pi}_{n+1}(Y)\ar[r]^{\Omega g_n}&\overline{\pi}_{n+1}(H)\ar[r]^{\partial_n}&\overline{\pi}_n(P_g)\ar[r]^{g^{'}_{n}} & \overline{\pi}_n(Y)  \ar[r]^{g_{n}}
                & \overline{\pi}^n(H)             }
$$ for $n\geq 0$.
\end{proof}
\end{proposition}

% Acknowledgements are put at the end of the introduction, in this way:
\section{Acknowledgements}
This work was supported by start-up research funds provided by BIMSA. The authors are grateful to BIMSA for its support and thank Stephen Theriault for his valuable suggestions and many help.

% \section{Second section}

% More math.

% If you are using hyperref, it will say
%   Package hyperref Warning: Token not allowed in a PDF string
% if you use math in a (sub)section title.  This is because it uses
% the title to produce bookmarks that can be shown in a side panel,
% and these can't use math.  You can use the \texampleorpdfstring command
% to give an alternate way to render these, possibly using unicode.
% \section{Section title with math \texorpdfstring{$\sum$}{}}

% More math.

\end{document}